\documentclass[moor]{informs1}            

\usepackage{natbib}
 \NatBibNumeric
 \def\bibfont{\small}%
 \bibpunct[, ]{[}{]}{,}{n}{}{,}%
 
 \usepackage[colorlinks=true,breaklinks=true,bookmarks=true,urlcolor=blue,  citecolor=blue,linkcolor=blue,bookmarksopen=false,draft=false]{hyperref}
 
\TheoremsNumberedBySection
\EquationsNumberedBySection 

\usepackage{amsmath,amstext,amsfonts,amssymb,mathrsfs,euscript}

\usepackage{algorithm,algorithmic}
\usepackage{graphicx,color,fullpage,setspace,subfigure,cases}

\TITLE{Empirical Dynamic Programming}
\ARTICLEAUTHORS{%
\AUTHOR{William B. Haskell}
\AFF{EE \& CS Departments, University of Southern California\\ \EMAIL{\tt wbhaskel@usc.edu}}
\AUTHOR{Rahul Jain\thanks{Corresponding author. This research was supported via NSF CAREER award CNS-0954116 and ONR Young Investigator Award N00014-12-1-0766.}}
\AFF{EE \& ISE Departments, University of Southern California\\ \EMAIL{\tt rahul.jain@usc.edu}}
\AUTHOR{Dileep  Kalathil}
\AFF{EE Department, University of Southern California\\ \EMAIL{\tt manisser@usc.edu}}
}

\begin{document}

\ABSTRACT{
We propose empirical dynamic programming algorithms for Markov decision processes (MDPs). In these algorithms, the exact expectation in the Bellman operator in classical value iteration is replaced by an empirical estimate to get `empirical value iteration' (EVI). Policy evaluation and policy improvement in classical policy iteration are also replaced by simulation to get `empirical policy iteration' (EPI). Thus, these empirical dynamic programming algorithms involve iteration of a random operator, the empirical Bellman operator. We introduce notions of probabilistic fixed points for such random monotone operators. We develop a stochastic dominance framework for convergence analysis of such operators. We then use this to give sample complexity bounds for both EVI and EPI. We then provide various variations and extensions to asynchronous empirical dynamic programming, the minimax empirical dynamic program, and show how this can also be used to solve the dynamic newsvendor problem. Preliminary experimental results suggest a faster rate of convergence than stochastic approximation algorithms.
}

\MSCCLASS{49L20, 90C39, 37M05, 62C12, 47B80, 37H99.}
\ORMSCLASS{TBD.}
\HISTORY{Submitted: \today}

\maketitle

\section{Introduction}\label{sec:intro}

Markov decision processes (MDPs) are natural models for decision making in a  stochastic dynamic setting for a wide variety of applications. The `principle of optimality' introduced by Richard Bellman in the 1950s has proved to be one of the most important ideas in stochastic control and optimization theory. It leads to dynamic programming algorithms for solving sequential stochastic optimization problems. And yet, it is well-known that it suffers from a ``curse of dimensionality'' \cite{Bellman57,nemhauser1966introduction,howard1971dynamic}, and does not scale computationally with state and action space size. In fact, the dynamic programming algorithm is known to be PSPACE-hard \cite{papadimitriou1987complexity}.

This realization led to the development of a wide variety of `approximate dynamic programming' methods beginning with the early work of Bellman himself \cite{bellman1959functional}. These ideas evolved independently in different fields, including the work of Werbos \cite{werbos1974beyond}, Kushner and Clark \cite{kushner1978stochastic} in control theory, Minsky \cite{minsky1961steps}, Barto, et al \cite{barto1981associative} and others in computer science, and Whitt in operations research \cite{whitt1978approximations,whitt1979approximations}. The key idea was an approximation of value functions using basis function approximation \cite{bellman1959functional}, state aggregation \cite{bertsekas1995dynamic},  and subsequently function approximation via neural networks \cite{anthony2009neural}. The difficulty was universality of the methods. Different classes of problems require different approximations.

Thus, alternative model-free methods were introduced by Watkins and Dayan \cite{watkins1992q} where a Q-learning algorithm was proposed as an approximation of the value iteration procedure. It was soon noticed that this is essentially a stochastic approximations scheme introduced in the 1950s by Robbins and Munro \cite{robbins1951stochastic}  and further developed by Kiefer and Wolfowitz \cite{kiefer1952stochastic} and Kushner and Clarke \cite{kushner1978stochastic}. This led to many subsequent generalizations including the temporal differences methods \cite{Bertsekas_Neuro_1996} and actor-critic methods \cite{Konda_Actor_1999, konda2004convergence}. These are summarized in  \cite{Bertsekas_Neuro_1996, sutton1998reinforcement, Powell_ADP_2007}. One shortcoming in this theory is that most of these algorithms require a recurrence property to hold, and in practice, often work only for finite state and action spaces. Furthermore, while many techniques for establishing convergence  have been developed \cite{borkar2008stochastic}, including the o.d.e. method \cite{ljung1977analysis, Borkar:2000}, establishing rate of convergence has been quite difficult \cite{Kakade_Reinforcement_2003}. Thus, despite considerable progress, these methods are not universal, sample complexity bounds are not known, and so other directions need to be explored.

A natural thing to consider is simulation-based methods. In fact, engineers and computer scientists often do dynamic programming via Monte-Carlo simulations. This technique affords considerable reduction in computation but at the expense of uncertainty about convergence. In this paper, we analyze such `empirical dynamic programming' algorithms. The idea behind the algorithms is quite simple and natural. In the DP algorithms, replace the expectation in the Bellman operator with a sample-average approximation. The idea is widely used in the stochastic programming literature, but mostly for single-stage problems. In our case, replacing the expectation with an empirical expectation operator, makes the classical Bellman operator a random operator. In the DP algorithm, we must find the fixed point of the Bellman operator. In the empirical DP algorithms, we must find a probabilistic fixed point of the random operator.

In this paper, we first introduce two notions of probabilistic fixed points that we call `strong' and `weak'. We then show that asymptotically
these concepts converge to the deterministic fixed point of the classical Bellman operator. The key technical idea of this paper is a novel stochastic dominance argument that is used to establish probabilistic convergence of a random operator, and in particular, of our empirical algorithms. Stochastic dominance, the notion of an order on the space of random variables, is a well developed tool (see \cite{muller2002comparison,shaked2007stochastic} for a comprehensive study).

In this paper, we develop a theory of empirical dynamic programming (EDP) for Markov decision processes (MDPs). Specifically, we make the following contributions in this paper. First, we propose both empirical value iteration and policy iteration algorithms and show that these converge. Each is an empirical variant of the classical algorithms. In empirical value iteration (EVI), the expectation in the Bellman operator is replaced with a sample-average or empirical approximation. In empirical policy iteration (EPI), both policy evaluation and the policy iteration are done via simulation, i.e., replacing the exact expectation with the simulation-derived empirical expectation. We note that the EDP class of algorithms is not a stochastic approximation scheme. Thus, we don't need a recurrence property as is commonly needed by stochastic approximation-based methods. Thus, the EDP algorithms are relevant for a larger class of problems (in fact, for any problem for which exact dynamic programming can be done.)  We provide convergence and sample complexity bounds for both EVI and EPI. But we note that EDP algorithms are essentially ``off-line'' algorithms just as classical DP is. Moreover, we also inherit some of the problems of classical DP such as scalability issues with large state spaces. These can be overcome in the same way as one does for classical DP, i.e., via state aggregation and function approximation.

Second, since the empirical Bellman operator is a random monotone operator, it doesn't have a deterministic fixed point. Thus, we introduce new mathematical notions of probabilistic fixed points. These concepts are pertinent when we are approximating a deterministic operator with an improving sequence of random operators. Under fairly mild assumptions, we show that our two probabilistic fixed point concepts converge to the deterministic fixed point of the classical monotone operator.

Third, since scant mathematical methods exist for convergence analysis of random operators, we develop a new technique based on stochastic dominance for convergence analysis of iteration of random operators. This technique allows for finite sample complexity bounds. We use this idea to prove convergence of the empirical Bellman operator by constructing a dominating Markov chain. We note that there is an extant theory of monotone random operators developed in the context of random dynamical systems \cite{chueshov2002monotone} but the techniques for convergence analysis of random operators is not relevant to our context. Our stochastic dominance argument can be applied for more general random monotone operators than just the empirical Bellman operator.

We also give a number of extensions of the EDP algorithms. We show that EVI can be performed asynchronously, making a parallel implementation possible. Second, we show that a saddle point equilibrium of a zero-sum stochastic game can be computed approximately by using the minimax Bellman operator. Third, we also show how the EDP algorithm and our convergence techniques can be used even with continuous state and action spaces by solving the dynamic newsvendor problem.

\subsubsection*{Related Literature}
A key question is how is the empirical dynamic programming method different from other methods for simulation-based optimization of MDPs, on which there is substantial literature. We note that most of these are stochastic approximation algorithms, also called reinforcement learning in computer science. Within this class. there are Q learning algorithms, actor-critic algorithms, and appproximate policy iteration algorithms. Q-learning was introduced by Watkins and Dayan but its convergence as a stochastic approximation scheme was done by Bertsekas and Tsitsiklis \cite{Bertsekas_Neuro_1996}. Q-learning for the average cost case was developed in \cite{abounadi2001learning}, and in a risk-sensitive context was developed in \cite{Borkar:2002uq}. The convergence rate of Q-learning was established in \cite{even2004learning} and similar sample complexity bounds were given in \cite{Kakade_Reinforcement_2003}. Actor-critic algorithms as a two time-scale stochastic approximation was developed in  \cite{Konda_Actor_1999}. But the most closely related work is optimistic policy iteration \cite{Tsitsiklis:2003} wherein simulated trajectories are used for policy evaluation while policy improvement is done exactly. The algorithm is a stochastic approximation scheme and its almost sure convergence follows. This is true of all stochastic approximation schemes but they do require some kind of recurrence property to hold. In contrast, EDP is not a stochastic approximation scheme, hence it does not need such assumptions. However, we can only guarantee its convergence in probability.

A class of simulation-based optimization algorithms for MDPs that is not based on stochastic approximations is the adaptive sampling methods developed by Fu, Marcus and co-authors \cite{chang2007simulation,chang2007survey}. These are based on the pursuit automata learning algorithms \cite{Thathachar85,rajaraman1996finite,narendra2012learning} and combine multi-armed bandit learning ideas with Monte-Carlo simulation to adaptively sample state-action pairs to approximate the value function of a finite-horizon MDP.

Some other closely related works are \cite{cooper2003convergence} which introduces simulation-based policy iteration (for average cost MDPs). It basically shows that almost sure convergence of such an algorithm can fail. Another related work is \cite{Cooper_Empirical_2011}. wherein a simulation-based value iteration is proposed for a finite horizon problem. Convergence in probability is established if the simulation functions corresponding to the MDP is Lipschitz continuous. Another closely related paper is \cite{Almudevar:2008fk}, which considers value iteration with error.   We note that our focus is on infinite horizon discounted MDPs. Moreover, we do not require any Lipschitz continuity condition. We show that EDP algorithms will converge (probabilistically) whenever the classical DP algorithms will converge (which is almost always).

A survey on approximate policy iteration is provided in \cite{Bertsekas_Policy_2011}. Approximate dynamic programming (ADP) methods are surveyed in \cite{Powell_ADP_2007}. In fact, many Monte-Carlo-based dynamic programming algorithms are introduced herein (but without convergence proof.) Simulation-based uniform estimation of value functions was studied in \cite{Jain_PAC_2006, Jain_PAC_2010}. This gave PAC learning type sample complexity bounds for MDPs and this can be combined with policy improvement along the lines of optimistic policy iteration.

This paper is organized as follows. In Section \ref{sec:prelim}, we discuss preliminaries and briefly talk about classical value and policy iteration. Section \ref{sec:edp} presents empirical value and policy iteration. Section \ref{sec:randomops} introduces the notion of random operators and relevant notions of probabilistic fixed points. In this section, we also develop a stochastic dominance argument for convergence analysis of iteration of random operators when they satisfy certain assumptions. In Section \ref{sec:samplecomplexity}, we show that the empirical Bellman operator satisfies the above assumptions, and present sample complexity and convergence rate estimates for the EDP algorithms. Section \ref{sec:variations} provides various extensions including asychronous EDP, minmax EDP and EDP for the dynamic newsvendor problem. Basic numerical experiments are reported in Section \ref{sec:experimental}.

\section{Preliminaries}\label{sec:prelim}

We first introduce  a typical representation
for a discrete time MDP as the 5-tuple
\[
\left(\mathbb{S},\,\mathbb{A},\,\left\{ A\left(s\right)\mbox{ : }s\in\mathbb{S}\right\} ,\, Q,\, c\right).
\]
Both the state space $\mathbb{S}$ and the action
space $\mathbb{A}$ are finite. Let $\mathcal{P}\left(\mathbb{S}\right)$ denote the space of probability
measures over $\mathbb{S}$ and we define $\mathcal{P}\left(\mathbb{A}\right)$ similarly.
For each state $s\in\mathbb{S}$, the set $A\left(s\right)\subset\mathbb{A}$
is the set of feasible actions. The entire set of feasible state-action
pairs is

\[
\mathbb{K}\triangleq\left\{ \left(s,a\right)\in\mathbb{S}\times\mathbb{A}\mbox{ : }a\in A\left(s\right)\right\} .
\]
The transition law $Q$ governs the system evolution, $Q\left(\cdot|s,a\right)\in\mathcal{P}\left(\mathbb{A}\right)$
for all $\left(s,\, a\right)\in\mathbb{K}$, i.e., $Q\left(j|s,a\right)$
for $j\in\mathbb{S}$ is the probability of visiting the state $j$
next given the current state-action pair $\left(s,\, a\right)$. Finally,
$c\mbox{ : }\mathbb{K}\rightarrow\mathbb{R}$ is a cost function that
depends on state-action pairs. 

Let $\Pi$ denote the class of \textit{stationary deterministic Markov policies}, i.e.,
mappings $\pi\mbox{ : }\mathbb{S}\rightarrow\mathbb{A}$ which only
depend on history through the current state. We only consider such
policies since it is well known that there is an optimal policy in
this class. For a given state $s\in\mathbb{S}$, $\pi\left(s\right)\in A\left(s\right)$
is the action chosen in state $s$ under the policy $\pi$. We assume
that $\Pi$ only contains feasible policies that respect the constraints
$\mathbb{K}$. The state and action at time $t$ are denoted $s_{t}$
and $a_{t}$, respectively. Any policy $\pi\in\Pi$ and initial state
$s\in\mathbb{S}$ determine a probability measure $P_{s}^{\pi}$ and
a stochastic process $\left\{ \left(s_{t},a_{t}\right),\, t\geq0\right\} $
defined on the canonical measurable space of trajectories of state-action
pairs. The expectation operator with respect to $P_{s}^{\pi}$ is
denoted $\mathbb{E}_{s}^{\pi}\left[\cdot\right]$.

We will focus on infinite horizon discounted cost MDPs with discount
factor $\alpha\in\left(0,\,1\right)$. For a given initial state $s\in\mathbb{S}$,
the expected discounted cost for policy $\pi\in\Pi$ is denoted by
\[
v^{\pi}\left(s\right)=\mathbb{E}_{s}^{\pi}\left[\sum_{t=0}^{\infty}\alpha^{t}c\left(s_{t},a_{t}\right)\right].
\]
The optimal cost starting from state $s$ is denoted by 
\[
v^{*}\left(s\right)\triangleq\inf_{\pi\in\Pi}\mathbb{E}_{s}^{\pi}\left[\sum_{t\geq0}\alpha^{t}c\left(s_{t},a_{t}\right)\right],
\]
and $v^{*}\in\mathbb{R}^{|\mathbb{S}|}$ denotes the corresponding
optimal value function in its entirety.

\paragraph*{Value iteration}
The Bellman operator $T\mbox{ : }\mathbb{R}^{|\mathbb{S}|}\rightarrow\mathbb{R}^{|\mathbb{S}|}$ is defined as
\[
\left[T\, v\right]\left(s\right)\triangleq\min_{a\in A\left(s\right)}\left\{ c\left(s,a\right)+\alpha\,\mathbb{E}\left[v\left(\tilde{s}\right)|s,a\right]\right\} ,\,\forall s\in\mathbb{S},
\]
for any $v\in\mathbb{R}^{|\mathbb{S}|}$, where $\tilde{s}$ is the
random next state visited, and
\[
\mathbb{E}\left[v\left(\tilde{s}\right)|s,a\right]=\sum_{j\in\mathbb{S}}v\left(j\right)Q\left(j|s,a\right)
\]
is the explicit computation of the expected cost-to-go conditioned
on state-action pair $\left(s,a\right)\in\mathbb{K}$. Value iteration
amounts to iteration of the Bellman operator. We have a sequence $\left\{ v^{k}\right\} _{k\geq0}\subset\mathbb{R}^{|\mathbb{S}|}$
where $v^{k+1}=T\, v^{k}=T^{k+1}v^{0}$ for all $k\geq0$ and an initial
seed $v^{0}$. This is the well-known value iteration algorithm for
dynamic programming.

We next state the Banach fixed point theorem which is used to prove
that value iteration converges. Let $U$ be a Banach space with norm
$\|\cdot\|_{U}$. We call an operator $G\mbox{ : }U\rightarrow U$
a $\textit{contraction mapping}$ when there exists a constant $\kappa\in[0,\,1)$
such that 
\begin{align*}
\|G\, v_{1}-G\, v_{2}\|_{U}\leq\kappa\,\|v_{1}-v_{2}\|_{U},\, & \forall v_{1},\, v_{2}\in U.
\end{align*}

\begin{theorem}
(Banach fixed point theorem) Let $U$ be a Banach space with norm
$\|\cdot\|_{U}$, and let $G\mbox{ : }U\rightarrow U$ be a contraction
mapping with constant $\kappa\in[0,\,1)$. Then,

(i) there exists a unique $v^{*}\in U$ such that $G\, v^{*}=v^{*}$;

(ii) for arbitrary $v^{0}\in U$, the sequence $v^{k}=G\, v^{k-1}=G^{k}v^{0}$
converges in norm to $v^{*}$ as $k\rightarrow\infty$;

(iii) $\|v^{k+1}-v^{*}\|_{U}\leq\kappa\,\|v^{k}-v^{*}\|_{U}$ for
all $k\geq0$.
\end{theorem}
For the rest of the paper, let $\mathcal{C}$ denote the space of
contraction mappings from $\mathbb{R}^{|\mathbb{S}|}\rightarrow\mathbb{R}^{|\mathbb{S}|}$.
It is well known that the Bellman operator $T\in\mathcal{C}$ with
constant $\kappa=\alpha$ is a contraction operator, and hence has
a unique fixed point $v^{*}$. It is known that value iteration converges
to $v^{*}$ as $k \rightarrow\infty$. In fact, $v^{*}$ is the optimal
value function. 


\paragraph*{Policy iteration}

Policy iteration is another well known dynamic programming algorithm
for solving MDPs. For a fixed policy $\pi\in\Pi$, define $T_{\pi}\mbox{ : }\mathbb{R}^{|\mathbb{S}|}\rightarrow\mathbb{R}^{|\mathbb{S}|}$
as
\[
\left[T_{\pi}v\right]\left(s\right)=c\left(s,\pi\left(s\right)\right)+\alpha\,\mathbb{E}\left[v\left(\tilde{s}\right)|s,\pi\left(s\right)\right].
\]
The first step is a $\textit{policy evalution}$ step. Compute $v^{\pi}$
by solving $T_{\pi}v^{\pi}=v^{\pi}$ for $v^{\pi}$. Let $c^{\pi}\in\mathbb{R}^{|\mathbb{S}|}$
be the vector of one period costs corresponding to a policy $\pi$,
$c^{\pi}\left(s\right)=c\left(s,\pi\left(s\right)\right)$ and $Q^{\pi}$,
the transition kernel corresponding to the policy $\pi$. Then, writing
$T_{\pi}v^{\pi}=v^{\pi}$ we have the linear system
\[
c^{\pi}+Q^{\pi}v^{\pi}=v^{\pi}.~~\text{(Policy Evaluation)}
\]
The second step is a $\textit{policy improvement}$ step. Given a
value function $v\in\mathbb{R}^{|\mathbb{S}|}$, find an `improved'
policy $\pi\in\Pi$ with respect to $v$ such that
\[
T_{\pi}v=T\, v. ~~\text{(Policy Update)}
\]
Thus, policy iteration produces a sequence of policies $\left\{ \pi^{k}\right\} _{k\geq0}$
and $\left\{ v^{k}\right\} _{k\geq0}$ as follows. At iteration $k \geq0$,
we solve the linear system $T_{\pi^{k}}v^{\pi^{k}}=v^{\pi^{k}}$ for
$v^{\pi^{k}}$, and then we choose a new policy $\pi^{k}$ satisfying

\[
T_{\pi^{k}}v^{\pi^{k}}=T\, v^{\pi^{k}},
\]
which is greedy with respect to $v^{\pi^{k}}$. We have a linear convergence
rate for policy iteration as well. Let $v\in\mathbb{R}^{|\mathbb{S}|}$
be any value function, solve $T_{\pi}v=T\, v$ for $\pi$, and then
compute $v^{\pi}$. Then, we know \cite[Lemma 6.2]{Bertsekas_Neuro_1996}
that
\[
\|v^{\pi}-v^{*}\|\leq\alpha\,\|v-v^{*}\|,
\]
from which convergence of policy iteration follows. Unless otherwise specified, the norm $||\cdot||$ we will use in this paper is the sup norm.

We use the following helpful fact in the paper. Proof is given in Appendix \ref{pf:fact}. 
\begin{remark}
\label{fact} Let $X$ be a given set, and $f_{1}\mbox{ : }X\rightarrow\mathbb{R}$
and $f_{2}\mbox{ : }X\rightarrow\mathbb{R}$ be two real-valued functions
on $X$. Then,

(i) $|\inf_{x\in X}f_{1}\left(x\right)-\inf_{x\in X}f_{2}\left(x\right)|\leq\sup_{x\in X}|f_{1}\left(x\right)-f_{2}\left(x\right)|$,
and

(ii) $|\sup_{x\in X}f_{1}\left(x\right)-\sup_{x\in X}f_{2}\left(x\right)|\leq\sup_{x\in X}|f_{1}\left(x\right)-f_{2}\left(x\right)|$.
\end{remark}

\section{Empirical Algorithms for Dynamic Programming}\label{sec:edp}

We now present empirical variants of dynamic programming algorithms. Our focus will be on value and policy iteration. As the reader will see, the idea is simple and natural. In subsequent sections we will introduce the new notions and techniques to prove their convergence.

\subsection{Empirical Value Iteration}\label{sec:evi}

We introduce   empirical value iteration (EVI) first. 

The Bellman operator $T$ requires exact evaluation of the expectation
\[
\mathbb{E}\left[v\left(\tilde{s}\right)|s,a\right]=\sum_{j\in\mathbb{S}}Q\left(j|s,a\right)v\left(j\right).
\]
We will simulate and replace this  exact expectation with an empirical estimate in each iteration. Thus, we need a simulation model for the MDP. Let
\[
\psi\mbox{ : }\mathbb{S}\times\mathbb{A}\times\left[0,1\right]\rightarrow\mathbb{S}
\]
be a simulation model for the state evolution for the MDP, i.e. $\psi$
yields the next state given the current state, the action taken and
an i.i.d. random variable. Without loss of generality, we can assume
that $\xi$ is a uniform random variable on $\left[0,1\right]$ and
$\left(s,a\right)\in\mathbb{K}$. With this convention, the Bellman
operator can be written as
\[
\left[T\, v\right]\left(s\right)\triangleq\min_{a\in A\left(s\right)}\left\{ c\left(s,a\right)+\alpha\,\mathbb{E}\left[v\left(\psi\left(s,a,\xi\right)\right)\right]\right\} ,\,\forall s\in\mathbb{S}.
\]

Now, we replace the expectation $\mathbb{E}\left[v\left(\psi\left(s,a,\xi\right)\right)\right]$
with its sample average approximation by simulating $\xi$. Given
$n$ i.i.d. samples of a uniform random variable, denoted $\left\{ \xi_{i}\right\} _{i=1}^{n}$,
the empirical estimate of $\mathbb{E}\left[v\left(\psi\left(s,a,\xi\right)\right)\right]$
is $\frac{1}{n}\sum_{i=1}^{n}v\left(\psi\left(s,a,\xi_{i}\right)\right)$.
We note that the samples are regenerated at each iteration. Thus,
the EVI algorithm can be summarized as follows.

\begin{algorithm}[!tph]
\caption{Empirical Value Iteration}

Input: $v^{0}\in\mathbb{R}^{|\mathbb{S}|}$, sample size $n\geq1$.

Set counter $k=0$.
\begin{enumerate}
\item Sample $n$ uniformly distributed random variables $\left\{ \xi_{i}\right\} _{i=1}^{n}$
from $\left[0,\,1\right]$, and compute
\[
v^{k+1}\left(s\right)=\min_{a\in A\left(s\right)}\left\{ c\left(s,a\right)+\frac{\alpha}{n}\sum_{i=1}^{n}v^{k}\left(\psi\left(s,a,\xi_{i}\right)\right)\right\} ,\,\forall s\in\mathbb{S}.
\]

\item Increment $k:=k+1$ and return to step 2.\end{enumerate}
\end{algorithm}


In each iteration, we regenerate samples and use this empirical estimate
to approximate $T$. Now we give  the sample complexity of the EVI algorithm. Proof is given in Section \ref{sec:samplecomplexity}.

\begin{theorem}
\label{thm:EVI_main} Given $\epsilon \in (0, 1)$ and $\delta \in (0, 1)$, fix $\epsilon_{g}=\epsilon/\eta^{*}$ and select $\delta_{1}, \delta_{2} > 0$ such that $\delta_{1} + 2 \delta_{2} \leq \delta$ where $\eta^{*}=\lceil 2/(1-\alpha) \rceil$. 
Select an $n$ such that
\[n \geq  n(\epsilon, \delta) = \frac{2\left(\kappa^{*}\right)^{2}}{\epsilon_{g}^{2}} \log \frac{ 2|\mathbb{K}|}{\delta_{1}}\]
where $\kappa^{*} = \max_{(s,a) \in \mathbb{K}} c(s,a)/(1-\alpha)$ and select a $k$ such that
\[ k \geq k(\epsilon, \delta)  = \log\left(\frac{1}{\delta_{2}\,\mu_{n,\,\min}}\right), \]
where $\mu_{n,min} = \min_{\eta} \mu_{n}(\eta)$ and   $\mu_{n}(\eta)$ is given by Lemma \ref{lem:steadystate}.
Then
\[
\mathcal{P}\left\{ \|\hat{v}_{n}^{k}-v^{*}\|\geq \epsilon\right\} \leq \delta.
\]
\end{theorem}

\begin{remark}
This result says that, if we take $n \geq n(\epsilon, \delta)$ samples in each iteration of the EVI algorithm and perform $k > k(\epsilon, \delta)$ iterations then the  EVI iterate $\hat{v}^{k}_{n}$ is $\epsilon$ close to the optimal value function $v^{*}$ with probability greater that $1-\delta$. We note that the sample complexity is $O\left(\frac{1}{\epsilon^{2}}, \log \frac{1}{\delta}, \log |\mathbb{S}|, \log |\mathbb{A}| \right)$. 
\end{remark}

The basic idea in the analysis is to frame EVI as iteration of a random operator $\widehat{T}_{n}$ which we call the empirical Bellman operator. We define $\widehat{T}_{n}$ as
\begin{equation}
\label{eq:empiricalBellman-defn}
\left[\widehat{T}_{n}\left(\omega\right)v\right]\left(s\right)\triangleq\min_{a\in A\left(s\right)}\left\{ c\left(s,a\right)+\frac{\alpha}{n}\sum_{i=1}^{n}v\left(\psi\left(s,a,\xi_{i}\right)\right)\right\} ,\,\forall s\in\mathbb{S}.
\end{equation}
This is a random operator because it depends on the random noise samples $\{\xi\}^{n}_{i=1}$. The definition and the analysis of this operator is done rigorously in Section \ref{sec:samplecomplexity}.

\subsection{Empirical Policy Iteration}\label{sec:epi}

We now define EPI along the same lines by replacing exact policy improvement
and evaluation with empirical estimates. For a fixed policy $\pi\in\Pi$,
we can estimate $v^{\pi}\left(s\right)$ via simulation. Given a sequence
of noise $\omega=\left(\xi_{i}\right)_{i\geq0}$, we have $s_{t+1}=\psi\left(s_{t},\pi\left(s_{t}\right),\xi_{t}\right)$
for all $t\geq0$. For $\gamma>0$, choose a finite horizon $\mathfrak{T}$
such that
\begin{equation*}
\label{eq:gamma-finitehorizon-defn}
\max_{\left(s,a\right)\in\mathbb{K}}|c\left(s,a\right)|\sum_{t=\mathfrak{T}+1}^{\infty}\alpha^{t}<\gamma.
\end{equation*}

We use the time horizon $\mathfrak{T}$ to truncate simulation, since
we must stop simulation after finite time. Let
\[
\left[\hat{v}^{\pi}\left(s\right)\right]\left(\boldsymbol{\omega}\right)=\sum_{t=0}^{\mathfrak{T}}\alpha^{t}c\left(s_{t}\left(\boldsymbol{\omega}\right),\pi\left(s_{t}\left(\boldsymbol{\omega}\right)\right)\right)
\]
be the realization of $\sum_{t=0}^{\mathfrak{T}}\alpha^{t}c\left(s_{t},a_{t}\right)$
on the sample path $\boldsymbol{\omega}$.

The next algorithm requires two input parameters, $n$ and $q$, which
determine sample sizes. Parameter $n$ is the sample size for policy
improvement and parameter $q$ is the sample size for policy evaluation.
We discuss the choices of these parameters in detail later. In the
following algorithm, the notation $s_{t}\left(\omega_{i}\right)$
is understood as the state at time $t$ in the simulated trajectory
$\omega_{i}$.

\begin{algorithm}[!tph]
\caption{Empirical Policy Iteration}

Input: $\pi_{0} \in \Pi$, $\epsilon>0$.
\begin{enumerate}
\item Set counter $k=0$.

\item For each $s\in\mathbb{S}$, draw $\omega_{1},\ldots,\omega_{q}\in\Omega$
and compute
\[
\hat{v}^{\pi_{k}}\left(s\right)=\frac{1}{q}\sum_{i=1}^{q}\sum_{t=0}^{\mathfrak{T}}\alpha^{t}c\left(s_{t}\left(\omega_{i}\right),\pi\left(s_{t}\left(\omega_{i}\right)\right)\right).
\]

\item Draw $\xi_{1},\ldots,\xi_{n}\in\left[0,1\right]$. Choose $\pi_{k+1}$
to satisfy
\[
\pi_{k+1}\left(s\right)\in\arg\min_{a\in A\left(s\right)}\left\{ c\left(s,a\right)+\frac{\alpha}{n}\sum_{i=1}^{n}\hat{v}^{\pi_{k}}\left(\psi\left(s,a,\xi_{i}\right)\right)\right\} ,\,\forall s\in\mathbb{S}.
\]

\item Increase $k:=k+1$ and return to step 2.
\item Stop when $\|\hat{v}^{\pi_{k+1}} - \hat{v}^{\pi_{k}} \|\leq\epsilon$.\end{enumerate}
\end{algorithm}

Step 2 replaces computation of $T_{\pi}v=T\, v$ (policy improvement).
Step 3 replaces solution of the system $v=c^{\pi}+\alpha\, Q^{\pi}v$
(policy evaluation). 

We now give the sample complexity result for EPI. Proof is given in Section \ref{sec:samplecomplexity}.
\begin{theorem}
\label{thm:EPI_main}
Given $\epsilon \in (0, 1)$, $\delta \in (0, 1)$ select  $\delta_{1}, \delta_{2} > 0$ such that $\delta_{1} + 2 \delta_{2} < \delta$. Also select   $\delta_{11}, \delta_{12} > 0$ such that $\delta_{11} + \delta_{12}  < \delta$. Then, select $\epsilon_{1}, \epsilon_{2} > 0$ such that $\epsilon_{g} = \frac{\epsilon_{2}+2\alpha \epsilon_{1}}{(1-\alpha)}$ where $\epsilon_{g}=\epsilon/\eta^{*},  \eta^{*}=\lceil 2/(1-\alpha) \rceil$. Then, select a $q$ and $n$ such that 
\[q \geq q(\epsilon, \delta) = \frac{2(\kappa^{*} (\mathfrak{T} +1))^{2}}{ (\epsilon_{1}-\gamma)^{2}} \log \frac{2 |\mathbb{S}|}{\delta_{11}}.  \]
\[n \geq n(\epsilon, \delta) = \frac{2\left(\kappa^{*}\right)^{2}}{(\epsilon_{2}/\alpha)^{2}} \log \frac{ 2|\mathbb{K}|}{\delta_{12}}.\]
where $\kappa^{*} = \max_{(s,a) \in \mathbb{K}} c(s,a)/(1-\alpha)$, and select a $k$ such that
\[ k   \geq k(\epsilon, \delta)= \log\left(\frac{1}{\delta_{2}\,\mu_{n,q,\,\min}}\right), \]
where $\mu_{n,q,min} = \min_{\eta} \mu_{n,q}(\eta)$ and $\mu_{n,q}(\eta)$ is given by equation \eqref{eq:mu-nq}.  Then,
\[
\mathcal{P}\left\{ \|{v}^{\pi_{k}}-v^{*}\|\geq \epsilon\right\} \leq \delta.
\]
\end{theorem}

\begin{remark}
This result says that, if we do $q \geq q(\epsilon, \delta)$ simulation runs for empirical policy evaluation, $n \geq n(\epsilon, \delta)$ samples for empirical policy update  and perform $k > k(\epsilon, \delta)$ iterations then the  true value $v^{\pi_{k}}$ of the policy $\pi_{k}$  will be $\epsilon$ close to the optimal value function $v^{*}$ with probability greater that $1-\delta$. We note that  $q$ is $O\left(\frac{1}{\epsilon^{2}}, \log \frac{1}{\delta},  \log |\mathbb{S}|\right)$ and  $n$  is $O\left(\frac{1}{\epsilon^{2}}, \log \frac{1}{\delta}, \log |\mathbb{S}|, \log |\mathbb{A}| \right)$. 
\end{remark}

\section{Iteration of Random Operators}
\label{sec:randomops}

The empirical Bellman operator $\widehat{T}_{n}$ we defined in equation \eqref{eq:empiricalBellman-defn}  is a random operator. When it operates on a vector, it yields a random vector.  When $\widehat{T}_{n}$ is iterated, it produces a stochastic process and we are interested in the possible convergence of this stochastic process. The underlying assumption is that the random operator $\widehat{T}_{n}$ is an `approximation' of a deterministic operator  $T$ such that $\widehat{T}_{n}$ converges to $T$ (in a sense we will shortly make precise) as $n$ increases. For example the empirical Bellman operator approximates the classical Bellman operator. We make this intuition mathematically rigorous in this section. The discussion in this section is not specific to the Bellman operator, but applies whenever a deterministic operator $T$ is being approximated by an improving sequence of random operators $\{ \widehat{T}_{n}\} _{n\geq1}$. 
%
%

\subsection{Probabilistic Fixed Points of Random Operators}
\label{sec:probfp}

In this subsection we formalize the definition of a random operator, denoted by $\widehat{T}_{n}$.

Since $\widehat{T}_{n}$ is a random operator, we need an appropriate
probability space upon which to define $\widehat{T}_{n}$. 
So, we define the sample space $\Omega=\left[0,1\right]^{\infty}$, the
$\sigma-$algebra $\mathcal{F}=\mathcal{B}^{\infty}$ where $\mathcal{B}$
is the inherited Borel $\sigma-$algebra on $\left[0,1\right]$, and
the probability distribution $P$ on $\Omega$ formed by an infinite
sequence of uniform random variables. The primitive uncertainties
on $\Omega$ are infinite sequences of uniform noise $\omega=\left(\xi_{i}\right)_{i\geq0}$
where each $\xi_{i}$ is an independent uniform random variable on
$\left[0,1\right]$. We view $\left(\Omega,\mathcal{F},\mathcal{P}\right)$
as the appropriate probability space on which to define iteration
of the random operators $\left\{ \widehat{T}_{n}\right\} _{n\geq1}$.

Next we define a composition of random operators, $\widehat{T}_{n}^{k}$, on the probability space $\left(\Omega^{\infty},\mathcal{F}^{\infty},\mathcal{P}\right)$,
for all $k\geq0$ and all $n\geq1$ where, 
\[
\widehat{T}_{n}^{k}\left(\boldsymbol{\omega}\right)v=\widehat{T}_{n}\left(\omega_{k-1}\right)\widehat{T}_{n}\left(\omega_{k-2}\right)\cdots\widehat{T}_{n}\left(\omega_{0}\right)v.
\]
Note that $\boldsymbol{\omega} \in \Omega^{\infty}$ is an infinite sequence $(\omega_{j})_{j \geq 0}$ where each $\omega_{j} = (\xi_{j,i})_{i \geq 0}$. Then we can define the iteration of $\widehat{T}_{n}$ with an initial seed $\hat{v}_{n}^{0}\in\mathbb{R}^{|\mathbb{S}|}$ (we use the hat notation to emphasize
that the iterates are random variables generated by the empirical operator) as
\begin{equation}
\label{eq:vhat-iter-defn}
\widehat{v}_{n}^{k+1}=\,  \widehat{T}_{n}\widehat{v}_{n}^{k} = \widehat{T}_{n}^{k} \widehat{v}_{n}^{0}
\end{equation}


Notice that we only iterate $k$ for fixed $n$. The sample size $n$
is constant in every stochastic process $\left\{ \hat{v}_{n}^{k}\right\} _{k\geq0}$,
where $\hat{v}_{n}^{k}=\widehat{T}_{n}^{k}\hat{v}^{0}$, for all $k\geq1$.
For a fixed $\hat{v}_{n}^{0}$, we can view all $\hat{v}_{n}^{k}$
as measurable mappings from $\Omega^{\infty}$ to $\mathbb{R}^{|\mathbb{S}|}$
via the mapping $\hat{v}_{n}^{k}\left(\boldsymbol{\omega}\right)=\widehat{T}_{n}^{k}\left(\boldsymbol{\omega}\right)\hat{v}_{n}^{0}$.

The relationship between the fixed points of the deterministic operator
$T$ and the probabilistic fixed points of the random operator $\left\{ \widehat{T}_{n}\right\} _{n\geq1}$
depends on how $\left\{ \widehat{T}_{n}\right\} _{n\geq1}$ approximates
$T$. Motivated by the relationship between the classical and the
empirical Bellman operator, we will make the following assumption. 

\begin{assumption}
\label{ass:fixedpoint} ~~$\mathcal{P}\left(\lim_{n\rightarrow\infty}\|\widehat{T}_{n}v-T\, v\|\geq\epsilon\right)=0$
 $\forall \epsilon>0$ and $\forall v\in\mathbb{R}^{\| \mathbb{S} \|}$. Also $T$ has a (possibly non-unique) fixed point $v^{*}$ such that $T v^{*} = v^{*}$.
\end{assumption}
Assumption \ref{ass:fixedpoint} is equivalent to $\lim_{n\rightarrow\infty}\widehat{T}_{n}\left(\omega\right)v=T\, v$
for $P-$almost all $\omega\in\Omega$. Here, we benefit from defining
all of the random operators $\left\{ \widehat{T}_{n}\right\} _{n\geq1}$
together on the sample space $\Omega=\left[0,1\right]^{\infty}$,
so that the above convergence statement makes sense.

\paragraph{Strong probabilistic fixed point:}

We now introduce a natural probabilistic fixed point notion for $\left\{ \widehat{T}_{n}\right\} _{n\geq1}$,
in analogy to the definition of a fixed point, $\|T\, v^{*}-v^{*}\|=0$
for a deterministic operator.
\begin{definition}
A vector $\hat{v} \in\mathbb{R}^{|\mathbb{S}|}$ is a strong probabilistic
fixed point for the sequence $\left\{ \widehat{T}_{n}\right\} _{n\geq1}$
if
\[
\lim_{n\rightarrow\infty}\mathcal{P}\left(\|\widehat{T}_{n}\hat{v}-\hat{v}\|>\epsilon\right)=0,~\forall \epsilon > 0. 
\]
\end{definition}
We note that the above notion is defined for a sequence of random
operators, rather than for a single random operator. 

\begin{remark}
\label{rem:epsdel-spfp}
We can give a slightly more general notion of a probabilistic fixed point which we call  $(\epsilon, \delta)$-strong probablistic fixed point. For a \emph{fixed} $(\epsilon, \delta)$, we say that a vector $\hat{v} \in\mathbb{R}^{|\mathbb{S}|}$ is  an  $(\epsilon, \delta)$-strong probabilistic fixed point if there exists an $n_{0}(\epsilon, \delta)$ such that for all $n \geq n_{0}(\epsilon, \delta)$ we get $\mathcal{P}\left(\|\widehat{T}_{n}\hat{v}-\hat{v}\|>\epsilon\right) < \delta$. Note that, all strong probabilistic fixed points satisfy this condition for arbitrary $(\epsilon, \delta)$ and hence are   $(\epsilon, \delta)$-strong probabilistic fixed points. However the converse need not be true. In many cases we may be looking for an $\epsilon$-optimal `solution' with a $1-\delta$  `probabilistic guarantee' where $(\epsilon, \delta)$ is fixed a priori. In fact, this would provide an approximation to the strong probabilistic fixed point of the sequence of operators.
\end{remark}

\paragraph{Weak probabilistic fixed point:}

It is well known that iteration of a deterministic contraction operator
converges to its fixed point. It is unclear whether a similar property
would hold for random operators, and whether they would converge to
the strong probabilistic fixed point of the sequence $\left\{ \widehat{T}_{n}\right\} _{n\geq1}$
in any way. Thus, we define an apparently weaker notion of a probabilistic
fixed point that explicitly considers iteration.
\begin{definition}
A vector $\hat{v}\in\mathbb{R}^{|\mathbb{S}|}$ is a weak probabilistic
fixed point for $\left\{ \widehat{T}_{n}\right\} _{n\geq1}$ if
\[
\lim_{n\rightarrow\infty}\limsup_{k\rightarrow\infty}\mathcal{P}\left(\|\widehat{T}_{n}^{k}v-\hat{v}\|>\epsilon\right)=0,~\forall \epsilon > 0,~\forall v \in \mathbb{R}^{|\mathbb{S}|}
\]
\end{definition}
We use $\limsup_{k\rightarrow\infty}\mathcal{P}\left(\|\widehat{T}_{n}^{k}v-\hat{v}\|>\epsilon\right)$
instead of $\lim_{k\rightarrow\infty}\mathcal{P}\left(\|\widehat{T}_{n}^{k}v-\hat{v}\|>\epsilon\right)$
because the latter limit may not exist for any fixed $n\geq1$. 

\begin{remark}
\label{rem:epsdel-wpfp}
Similar to the definition that we gave in Remark \ref{rem:epsdel-spfp}, we can define an $(\epsilon, \delta)$-weak probablistic fixed point.  For a \emph{fixed} $(\epsilon, \delta)$, we say that a vector $\hat{v} \in\mathbb{R}^{|\mathbb{S}|}$ is  an  $(\epsilon, \delta)$-weak probabilistic fixed point if there exists an $n_{0}(\epsilon, \delta)$ such that for all $n \geq n_{0}(\epsilon, \delta)$ we get $\limsup_{k\rightarrow\infty}\mathcal{P}\left(\|\widehat{T}_{n}^{k}v-\hat{v}\|>\epsilon\right) < \delta$. As before, all weak probabilistic fixed points are indeed $(\epsilon, \delta)$ weak probabilistic fixed points, but converse need not be true. 
\end{remark}

At this point the connection between strong/weak probabilistic fixed points of the random operator $\widehat{T}_{n}$ and  the classical fixed point of the deterministic operator $T$ is not clear. Also it is not clear  whether the random sequence  $\{\hat{v}^{k}_{n}\}_{k \geq 0}$ converges to either of these two fixed point. In the following subsections we address these issues.

\subsection{A Stochastic Process on $\mathbb{N}$}\label{sec:stochproc}

In this subsection, we construct a new stochastic process on $\mathbb{N}$
that will be useful in our analysis.  We first start with a simple lemma.
\begin{lemma}
The stochastic process $\left\{ \hat{v}_{n}^{k}\right\}_{k\geq0}$
is a Markov chain on $\mathbb{R}^{|\mathbb{S}|}$.
\end{lemma}

\proof{Proof:}
This follows from the fact that each iteration of $\widehat{T}_{n}$ is
independent, and identically distributed. Thus, the next iterate $\hat{v}_{n}^{k+1}$
only depends on history through the current iterate $\hat{v}_{n}^{k}$. 
\Halmos
\endproof
\paragraph{}
Even though $\left\{ \hat{v}_{n}^{k}\right\}_{k\geq0}$ is a Markov
chain, its analysis is complicated by two factors. First, $\left\{ \hat{v}_{n}^{k}\right\} _{k\geq0}$
is a Markov chain on the continuous state space $\mathbb{R}^{|\mathbb{S}|}$,
which introduces technical difficulties in general when compared
to a discrete state space. Second, the transition probabilities of
$\left\{ \hat{v}_{n}^{k}\right\} _{k\geq0}$ are too complicated to
compute explicitly. 

Since we are approximating $T$ by $\widehat{T}_{n}$  and want
to compute $v^{*}$, we should track the progress of $\left\{ \hat{v}_{n}^{k}\right\} _{k\geq0}$
to the fixed point $v^{*}$ of $T$. Equivalently, we are interested
in the real-valued stochastic process $\left\{ \|\hat{v}_{n}^{k}-v^{*}\|\right\} _{k\geq0}$. 
If $\|\hat{v}_{n}^{k}-v^{*}\|$ approaches zero then $\hat{v}_{n}^{k}$
approaches $v^{*}$, and vice versa.

The state space of the stochastic process $\left\{ \|\hat{v}_{n}^{k}-v^{*}\|\right\} _{k\geq0}$
is $\mathbb{R}$, which is simpler than the state space $\mathbb{R}^{|\mathbb{S}|}$
of $\left\{ \hat{v}_{n}^{k}\right\} _{k\geq0}$, but which is 
still continuous. Moreover, $\left\{ \|\hat{v}_{n}^{k}-v^{*}\|\right\} _{k\geq0}$
is a non-Markovian process in general. In fact it would be easier to  study a related stochastic process on a discrete, and  ideally a finite state space. In this subsection we show how this can be done. 

We make a boundedness assumption next.
\begin{assumption}
\label{ass:bddness}
There exists a $\kappa^{*} < \infty$ such that $\|\hat{v}_{n}^{k}\| \leq  \kappa^{*}$ almost surely  for all $k \geq 0$, $n \geq 1$. Also, $\|v^{*}\| \leq \kappa^{*}$. 
\end{assumption}

Under this assumption we can  restrict the stochastic process $\left\{ \|\hat{v}^{k}-v^{*}\|\right\} _{k\geq0}$
to the compact state space
\[
\overline{B_{2\kappa^{*}}\left(0\right)}=\left\{ v\in\mathbb{R}^{|\mathbb{S}|}\mbox{ : }\|v\|\leq 2\kappa^{*}\right\} .
\]
We will adopt the convention that
any element $v$ outside of $\overline{B_{\kappa^{*}}\left(0\right)}$
will be mapped to its projection $\kappa^{*}\frac{v}{\|v\|}$ onto
$\overline{B_{\kappa^{*}}\left(0\right)}$ by any realization of $\widehat{T}_{n}$.

Choose a \textbf{granularity} $\epsilon_{g}>0$ to be fixed for the
remainder of this discussion. We will break up $\mathbb{R}$ into
intervals of length $\epsilon_{g}$ starting at zero, and we will
note which interval is occupied by $\|\hat{v}_{n}^{k}-v^{*}\|$ at
each $k\geq0$. We will define a new stochastic process $\left\{ X_{n}^{k}\right\} _{k\geq0}$
on $\left(\Omega^{\infty},\mathcal{F}^{\infty},\mathcal{P}\right)$
with state space $\mathbb{N}$. The idea is that $\left\{ X_{n}^{k}\right\} _{k\geq0}$
will report which interval of $\left[0,\,2\,\kappa^{*}\right]$ is
occupied by $\left\{ \|\hat{v}_{n}^{k}-v^{*}\|\right\} _{k\geq0}$.
Define $X_{n}^{k}\mbox{ : }\Omega^{\infty}\rightarrow\mathbb{N}$
via the rule:
\begin{equation}
X_{n}^{k}\left(\boldsymbol{\omega}\right)=\begin{cases}
0, & \mbox{if }\|\hat{v}^{k}\left(\boldsymbol{\omega}\right)-v^{*}\|=0,\\
\eta\geq1, & \mbox{if }\left(\eta-1\right)\epsilon_{g}<\|\hat{v}^{k}\left(\boldsymbol{\omega}\right)-v^{*}\|\leq\eta\,\epsilon_{g},
\end{cases}\label{eq:Xt}
\end{equation}
for all $k\geq0$. More compactly, 
\[
X_{n}^{k}\left(\boldsymbol{\omega}\right)=\left\lceil \|\hat{v}_{n}^{k}\left(\boldsymbol{\omega}\right)-v^{*}\|/\epsilon_{g}\right\rceil ,
\]
where $\left\lceil \chi\right\rceil $ denotes the smallest integer
greater than or equal to $\chi\in\mathbb{R}$. Thus  the stochastic process $\left\{ X_{n}^{k}\right\} _{k\geq0}$ is a
report on how close the stochastic process $\left\{ \|\hat{v}_{n}^{k}-v^{*}\|\right\} _{k\geq0}$
is to zero, and in turn how close the Markov chain $\left\{ \hat{v}_{n}^{k}\right\} _{k\geq0}$
is to the true fixed point $v^{*}$ of $T$.

Define the constant
\[
N^{*}\triangleq\left\lceil \frac{2\,\kappa^{*}}{\epsilon_{g}}\right\rceil .
\]
Notice $N^{*}$ is the smallest number of intervals of length $\epsilon_{g}$
needed to cover the interval $\left[0,\,2\,\kappa^{*}\right]$. By
construction, the stochastic process $\left\{ X_{n}^{k}\right\} _{k\geq0}$
is restricted to the finite state space $\left\{ \eta\in\mathbb{N}\mbox{ : }0\leq\eta\leq N^{*}\right\}.$

The process $\left\{ X_{n}^{k}\right\} _{k\geq0}$ need not be a
Markov chain. However, it is easier to work with than either $\left\{ \hat{v}_{n}^{k}\right\} _{k\geq0}$
or $\left\{ \|\hat{v}_{n}^{k}-v^{*}\|\right\} _{k\geq0}$ because
it has a discrete state space. It is also easy to relate $\left\{ X_{n}^{k}\right\} _{k\geq0}$
back to $\left\{ \|\hat{v}_{n}^{k}-v^{*}\|\right\} _{k\geq0}$. 

Recall that $X\geq_{as}Y$ denotes almost sure inequality between
two random variables $X$ and $Y$ defined on the same probability
space. The stochastic processes $\left\{ X_{n}^{k}\right\} _{k\geq0}$
and $\left\{ \|\hat{v}_{n}^{k}-v^{*}\|/\epsilon_{g}\right\} _{k\geq0}$
are defined on the same probability space, so the next lemma follows
by construction of $\left\{ X_{n}^{k}\right\} _{k\geq0}$.

\begin{lemma}
For all $k\geq0$, $X_{n}^{k}\geq_{as}\|\hat{v}_{n}^{k}-v^{*}\|/\epsilon_{g}$.
\end{lemma}

To proceed, we will  make the following assumptions  about the  deterministic operator $T$ and the random operator $\widehat{T}_{n}$. 
\begin{assumption}
\label{ass:rate}  $\|T\, v-v^{*}\|\leq\alpha\,\|v-v^{*}\|$
for all $v\in\mathbb{R}^{|\mathbb{S}|}$.
\end{assumption}

\begin{assumption}
\label{ass:pn}  There is a sequence $\{p_{n}\}_{n \geq 1}$ such that 
\[P\left(\|T v - \widehat{T}_{n}v \| < \epsilon \right) > p_{n}(\epsilon)  \]
and $p_{n}(\epsilon) \uparrow 1$ as $n \rightarrow \infty$ for all $v\in \overline{B_{\kappa^{*}}(0)}$, $\forall \epsilon > 0$ . 
\end{assumption}


We now discuss the convergence rate of $\left\{ X_{n}^{k}\right\} _{k\geq0}$. Let $X_{n}^{k}=\eta$. On the event $F=\left\{ \|T\,\hat{v}_{n}^{k}-\widehat{T}_{n}\hat{v}_{n}^{k}\|<\epsilon_{g}\right\}$, we have
\[\|\hat{v}_{n}^{k+1}-v^{*}\|\leq\,  \|\widehat{T}_{n}\hat{v}_{n}^{k}-T\,\hat{v}^{k}\| +  \|T\,\hat{v}_{n}^{k}-v^{*}\|
\leq\, \left(\alpha\,\eta+1\right)\epsilon_{g}\]
where we used Assumption \ref{ass:rate} and the definition of $ X_{n}^{k}$. Now using Assumption \ref{ass:pn} we can summarize:
\begin{equation}
\label{eq:xk-to-yk-step-evi}
\text{If}~X_{n}^{k}=\eta,~\text{then}~X_{n}^{k+1} \leq \left\lceil \alpha\,\eta+1\right\rceil ~\text{with a probability at least}~p_{n}(\epsilon_{g}).\
\end{equation}

We conclude this subsection with a comment about the state space of the stochastic process of $\{X_{n}^{k}\}_{k \geq 0}$.  If we start with $X_{n}^{k}=\eta$ and  if $\left\lceil \alpha\,\eta+1\right\rceil < \eta$ then
we must have improvement in the proximity of $\hat{v}_{n}^{k+1}$
to $v^{*}$. We define a new constant
\[
\eta^{*}=\min\left\{ \eta\in\mathbb{N}\mbox{ : }\left\lceil \alpha\,\eta+1\right\rceil < \eta\right\} =\left\lceil \frac{2}{1-\alpha}\right\rceil .
\]
If $\eta$ is too small, then $\left\lceil \alpha\,\eta+1\right\rceil $
may be equal to $\eta$ and no improvement in the proximity of $\hat{v}_{n}^{k}$
to $v^{*}$ can be detected by $\left\{ X_{n}^{k}\right\} _{k\geq0}$.
For any $\eta\geq\eta^{*}$, $\left\lceil \alpha\,\eta+1\right\rceil <\eta$
and strict improvement must hold. So, for the stochastic process $\left\{ X_{n}^{k}\right\} _{k\geq0}$, we can restrict our attention to the state space $\mathcal{X}:= \{\eta^{*}, \eta^{*}+1, \ldots, N^{*}-1, N^{*}\}$.

\subsection{Dominating Markov Chains}\label{sec:dominatingMC}

If we could understand the behavior of the stochastic processes $\left\{ X_{n}^{k}\right\} _{k\geq0}$,
then we could make statements about the convergence of $\left\{ \|\hat{v}_{n}^{k}-v^{*}\|\right\} _{k\geq0}$
and $\left\{ \hat{v}_{n}^{k}\right\} _{k\geq0}$. Although  simpler
than $\left\{ \|\hat{v}_{n}^{k}-v^{*}\|\right\} _{k\geq0}$ and $\left\{ \hat{v}_{n}^{k}\right\} _{k\geq0}$,
the stochastic process $\left\{ X_{n}^{k}\right\} _{k\geq0}$ is still
too complicated to work with analytically. We overcome this difficulty
with a family of dominating Markov chains.  We now present our dominance argument. Several technical
details are expanded upon in the appendix.

We will denote our family of ``dominating'' Markov chains (MC) by $\left\{ Y_{n}^{k}\right\} _{k\geq0}$.
We will construct these Markov chains to be tractable and to help
us analyze $\left\{ X_{n}^{k}\right\} _{k\geq0}$. Notice that the
family $\left\{ Y_{n}^{k}\right\} _{k\geq0}$ has explicit dependence
on $n\geq1$. We do not necessarily construct $\left\{ Y_{n}^{k}\right\} _{k\geq0}$
on the probability space $\left(\Omega^{\infty},\mathcal{F}^{\infty},\mathcal{P}\right)$.
Rather, we view $\left\{ Y_{n}^{k}\right\} _{k\geq0}$ as being defined
on $\left(\mathbb{N}^{\infty},\,\mathcal{N}\right)$, the canonical
measurable space of trajectories on $\mathbb{N}$, so $Y_{n}^{k}\mbox{ : }\mathbb{N}^{\infty}\rightarrow\mathbb{N}$.
We will use $\mathcal{Q}$ to denote the probability measure of $\left\{ Y_{n}^{k}\right\} _{k\geq0}$
on $\left(\mathbb{N}^{\infty},\,\mathcal{N}\right)$. Since $\left\{ Y_{n}^{k}\right\} _{k\geq0}$
will be a Markov chain by construction, the probability measure $\mathcal{Q}$
will be completely determined by an initial distribution on $\mathbb{N}$
and a transition kernel. We denote the transition kernel of $\left\{ Y_{n}^{k}\right\} _{k\geq0}$
as $\mathfrak{Q}_{n}$.

Our specific choice for $\left\{ Y_{n}^{k}\right\} _{k\geq0}$ is
motivated by analytical expediency, though the reader will see that many
other choices are possible. We now construct the process $\left\{ Y_{n}^{k}\right\} _{k\geq0}$
explicitly, and then  compute its steady state probabilities and
its mixing time. We will define the stochastic process $\left\{ Y_{n}^{k}\right\} _{k\geq0}$ on the finite
state space $\mathcal{X}$, 
based on our observations about the boundedness of $\left\{ \|\hat{v}_{n}^{k}-v^{*}\|\right\} _{k\geq0}$
and $\left\{ X_{n}^{k}\right\} _{k\geq0}$. Now, for a fixed $n$ and $p_{n}(\epsilon_{g})$ (we drop the argument $\epsilon_{g}$ in the following for notational convenience) as assumed in Assumption \ref{ass:pn} we construct the dominating Markov chain $\left\{ Y_{n}^{k}\right\} _{k\geq0}$ as:
\begin{equation}
\label{eq:Ynk-evi-construction }
Y_{n}^{k+1}=\begin{cases}
\max\left\{ Y_{n}^{k}-1,\,\eta^{*}\right\} , & \mbox{w.p. }p_{n},\\
N^{*}, & \mbox{w.p. }1-p_{n}.
\end{cases}
\end{equation}

The first value $Y_{n}^{k+1}=\max\left\{ Y_{n}^{k}-1,\,\eta^{*}\right\} $
corresponds to the case where the approximation error satisfies $\|T\,\hat{v}^{k}-\widehat{T}_{n}\hat{v}^{k}\|<\epsilon_{g}$,
and the second value $Y_{n}^{k+1}=N^{*}$ corresponds to all other
cases (giving us an extremely conservative bound in the sequel). This
construction also ensures that $Y_{n}^{k} \in \mathcal{X}$
for all $k$, $\mathcal{Q}-$almost surely. Informally, $\left\{ Y_{n}^{k}\right\} _{k\geq0}$
will either move one unit closer to zero until it reaches $\eta^{*}$,
or it will move (as far away from zero as possible) to $N^{*}$.

We now summarize some key properties of $\left\{ Y_{n}^{k}\right\} _{k\geq0}$.
\begin{proposition}
\label{prop:steadystateY}
For $\left\{ Y_{n}^{k}\right\} _{k\geq0}$ as defined above,

(i) it is a Markov chain;

(ii) the steady state distribution of $\left\{ Y_{n}^{k}\right\} _{k\geq0}$,
and the  limit $Y_{n}=_{d}\lim_{k\rightarrow\infty}Y^{k}_{n}$, exists;

(iii) $\mathcal{Q}\left\{ Y_{n}^{k} > \eta\right\} \rightarrow\mathcal{Q}\left\{ Y_{n} > \eta\right\} $
as $k\rightarrow\infty$ for all $\eta\in\mathbb{N}$;


\end{proposition}
\proof{Proof:}
Parts (i) - (iii) follow by construction of $\left\{ Y_{n}^{k}\right\} _{k\geq0}$
and the fact that this family consists of irreducible Markov chains
on a finite state space. 
\hfill \Halmos \endproof 


We now describe a stochastic dominance relationship between the two
stochastic processes $\left\{ X_{n}^{k}\right\} _{k\geq0}$ and $\left\{ Y_{n}^{k}\right\} _{k\geq0}$.
The notion of stochastic dominance (in the usual sense) will be central
to our development.
\begin{definition}
Let $X$ and $Y$ be two real-valued random variables, then $Y$ \textit{stochastically dominates}
$X$, written $X\leq_{st}Y$, when $\mathbb{E}\left[f\left(X\right)\right]\leq\mathbb{E}\left[f\left(Y\right)\right]$
for all increasing functions $f\mbox{ : }\mathbb{R}\rightarrow\mathbb{R}$.
\end{definition}
The condition $X\leq_{st}Y$ is known to be equivalent to
\[
\mathbb{E}\left[\boldsymbol{1}\left\{ X\geq\theta\right\} \right]\leq\mathbb{E}\left[\boldsymbol{1}\left\{ Y\geq\theta\right\} \right]~~\text{or}~~\mbox{Pr}\left\{ X\geq\theta\right\} \leq\mbox{Pr}\left\{ Y\geq\theta\right\}, 
\]
for all $\theta$ in the support of $Y$. Notice that the relation
$X\leq_{st}Y$ makes no mention of the respective probability spaces
on which $X$ and $Y$ are defined - these spaces may be the same
or different (in our case they are different).

Let $\left\{ \mathcal{F}^{k}\right\} _{k\geq0}$ be the filtration
on $\left(\Omega^{\infty},\mathcal{F}^{\infty},\mathcal{P}\right)$
corresponding to the evolution of information about $\left\{ X_{n}^{k}\right\} _{k\geq0}$.
Let $\left[X_{n}^{k+1}|\mathcal{F}^{k}\right]$ denote the
conditional distribution of $X_{n}^{k+1}$ given the information $\mathcal{F}^{k}$.
The following theorem compares the marginal distributions of $\left\{ X_{n}^{k}\right\} _{k\geq0}$
and $\left\{ Y_{n}^{k}\right\} _{k\geq0}$ at all times $k\geq0$
when the two stochastic processes $\left\{ X_{n}^{k}\right\} _{k\geq0}$
and $\left\{ Y_{n}^{k}\right\} _{k\geq0}$ start from the same state.

\begin{theorem}
\label{thm:dominance} 
If $X_{n}^{0}=Y_{n}^{0}$, then $X_{n}^{k}\leq_{st}Y_{n}^{k}$
for all $k\geq0$.
\end{theorem}

Proof is given in Appendix \ref{pf:thm:dominance}

%

The following corollary resulting from Theorem \ref{thm:dominance}
relates the stochastic processes $\left\{ \|\hat{v}_{n}^{k}-v^{*}\|\right\} _{k\geq0}$,
$\left\{ X_{n}^{k}\right\} _{k\geq0}$, and $\left\{ Y_{n}^{k}\right\} _{k\geq0}$
in a probabilistic sense, and summarizes our stochastic dominance
argument.

\begin{corollary}
\label{cor:probabilistic} For any fixed $n \geq 1$,   we have\\
(i)$\mathcal{P}\left\{ \|\hat{v}_{n}^{k}-v^{*}\| > \eta \epsilon_{g}\right\} \leq\mathcal{P}\left\{ X_{n}^{k} > \eta\right\} \leq\mathcal{Q}\left\{ Y_{n}^{k} > \eta\right\} $
for all $\eta\in\mathbb{N}$ for all $k \geq0$;\\
(ii) $\limsup_{k\rightarrow\infty}\mathcal{P}\left\{ X_{n}^{k} > \eta\right\} \leq\mathcal{Q}\left\{ Y_{n} > \eta\right\} $ for all $\eta\in\mathbb{N}$;\\
(iii) \textup{$\limsup_{k\rightarrow\infty}\mathcal{P}\left\{ \|\hat{v}_{n}^{k}-v^{*}\| > \eta\,\epsilon_{g}\right\} \leq\mathcal{Q}\left\{ Y_{n} > \eta\right\}$
for all }$\eta\in\mathbb{N}$.
\end{corollary}
\proof{Proof:}
(i) The first inequality is true by construction of $X^{k}_{n}$. Then $\mathcal{P}\left\{ X_{n}^{k} > \eta\right\} \leq\mathcal{Q}\left\{ Y_{n}^{k} > \eta\right\} $
for all $k \geq0$ and $\eta\in\mathbb{N}$ by Theorem \ref{thm:dominance}.\\
(ii) Since $\mathcal{Q}\left\{ Y_{n}^{k} > \eta\right\} $ converges (by Proposition \ref{prop:steadystateY}), the result follows by taking the limit in part (i).\\
(iii) This again follows by taking limit in part (i) and using Proposition \ref{prop:steadystateY}
\hfill \Halmos \endproof

We now compute the steady  state distribution of the Markov chain $\{Y_{n}^{k}\}_{k \geq 0}$. Let $\mu_{n}$ denotes the steady state distribution of  $Y_{n}=_{d}\lim_{k\rightarrow\infty}Y^{k}_{n}$ (whose existence is guaranteed by  Proposition \ref{prop:steadystateY}) where $\mu_{n}\left(i\right)=\mathcal{Q}\left\{ Y_{n}=i\right\} $
for all $i \in \mathcal{X}$.  The next lemma follows from standard techniques (see \cite{Ross_StochasticProcesses_1996} for example). Proof is given in  Appendix \ref{pf:lem:steadystate}.
\begin{lemma}
\label{lem:steadystate}
For any fixed $n \geq 1$, 
\[\mu_{n}\left(\eta^{*}\right)=\,~p_{n}^{N^{*}-\eta^{*}-1}, \mu_{n}\left(N^{*}\right)=\,  \frac{1-p_{n}}{p_{n}},~ \mu_{n}\left(i\right)=\,  \left(1-p_{n}\right)p_{n}^{\left(N^{*}-i-1\right)},  \forall i=\eta^{*}+1,\ldots,N^{*}-1 .  \]
\end{lemma}
Note that an explicit expression for $p_{n}$ in the case of empirical Bellman operator is given in equation \eqref{eq:defn-pn}.

\subsection{Convergence Analysis of Random Operators}\label{sec:convergence}

We now give  results on the convergence of the stochastic process $\{\hat{v}^{k}_{n}$, which could equivalently be written $\widehat{T}^{k}_{n}\hat{v}_{0}\}_{k \geq 0}$. Also we elaborate on the connections between our different notions of fixed points. Throughout this section, $v^{*}$ denotes the fixed point of the deterministic operator $T$ as defined in Assumption \ref{ass:fixedpoint}.

\begin{theorem}
\label{thm:asymptotic} Suppose  the  random operator $\widehat{T}_{n}$ satisfies the assumptions \ref{ass:fixedpoint} - \ref{ass:pn}. Then for any $\epsilon>0$,
\[
\lim_{n\rightarrow\infty}\limsup_{k\rightarrow\infty}\mathcal{P}\left(\|\hat{v}_{n}^{k}-v^{*}\|>\epsilon\right)=0,
\]
i.e. $v^{*}$ is a weak probabilistic fixed point of $\left\{ \widehat{T}_{n}\right\} _{n\geq1}$.
\end{theorem}
\proof{Proof:}
Choose  the granularity $\epsilon_{g} = \epsilon/\eta^{*}$. From Corollary \ref{cor:probabilistic} and Lemma \ref{lem:steadystate},
\[\limsup_{k\rightarrow\infty}\mathcal{P}\left\{ \|\hat{v}_{n}^{k}-v^{*}\|> \eta^{*}\epsilon_{g}\right\}  \leq\mathcal{Q}\left\{ Y_{n} > \eta^{*}\right\} = 1 - \mu_{n}\left(\eta^{*}\right) = 1 - p_{n}^{N^{*}-\eta^{*}-1} \]
Now by Assumption \ref{ass:pn}, $p_{n} \uparrow 1$ and by taking limit on both sides of the above inequality gives the desired result.  
\hfill \Halmos \endproof

Now we show that  a strong probabilistic fixed point and the deterministic fixed point $v^{*}$ coincide under Assumption \ref{ass:fixedpoint}. 

\begin{proposition}
\label{lem:strongpfp}
Suppose Assumption\textbf{ \ref{ass:fixedpoint}} holds. Then,  \\
(i)  $v^{*}$ is a strong probabilistic
fixed point of the sequence $\left\{ \widehat{T}_{n}\right\} _{n\geq1}$.\\
(ii) Let $\hat{v}$ be a strong probabilistic
fixed point of the sequence $\left\{ \widehat{T}_{n}\right\} _{n\geq1}$,
then $\hat{v}$ is a fixed point of $T$.
\end{proposition}
Proof is given in Appendix \ref{pf:lem:strongpfp}.

Thus the set of fixed points of $T$ and the set of strong probabilistic fixed points of $\{\widehat{T}_{n}\}_{n \geq 1}$ coincide. This suggests that a ``probabilistic'' fixed point would be an ``approximate'' fixed point of the deterministic operator $T$.

We now explore the connection between weak
probabilistic fixed points and classical fixed points.

\begin{proposition}
\label{prop:wkpfp1}Suppose  the  random operator $\widehat{T}_{n}$ satisfies the assumptions \ref{ass:fixedpoint} - \ref{ass:pn}. Then,\\
(i)  $v^{*}$ is a weak probabilistic fixed point of the
sequence $\left\{ \widehat{T}_{n}\right\} _{n\geq1}$. \\
(ii) Let $\hat{v}$ be a weak probabilistic fixed point of the sequence $\left\{ \widehat{T}_{n}\right\} _{n\geq1}$,
then $\hat{v}$ is a fixed point of $T$.
\end{proposition}
Proof is given Appendix \ref{pf:prop:wkpfp1}. It is obvious that we need more assumptions to analyze the asymptotic behavior of the iterates of the random operator $\widehat{T}_{n}$ and establish the connection to  the fixed point of the deterministic operator.

We summarize the above discussion in the following theorem.
\begin{theorem}
Suppose   the  random operator $\widehat{T}_{n}$ satisfies the assumptions \ref{ass:fixedpoint} - \ref{ass:pn}. Then the following three statements are equivalent:

(i) $v$ is a fixed point of $T$,

(ii) $v$ is a strong probabilistic fixed point of $\left\{ \widehat{T}_{n}\right\} _{n\geq1}$,

(iii) $v$ is a weak probabilistic fixed point of $\left\{ \widehat{T}_{n}\right\} _{n\geq1}$.
\end{theorem}

This is quite remarkable because we see not only that the two notions of a probabilistic fixed point of a sequence of random operators coincide, but in fact they coincide with the fixed point of the related classical operator. Actually, it would have been disappointing if this were not the case.  The above result now suggests that the iteration of a random operator a finite number $k$ of times and for a fixed $n$ would yield an approximation to the classical fixed point with high probability. Thus, the notions of  the $(\epsilon,\delta)$-strong and weak probabilistic fixed points coincide asymptotically, however, we note that non-asymptotically they need not be the same.

\section{Sample Complexity for EDP}\label{sec:samplecomplexity}

In this section we present the proofs of the sample complexity results
for empirical value iteration (EVI) and policy iteration (EPI) (Theorem \ref{thm:EVI_main} and Theorem \ref{thm:EPI_main} in Section \ref{sec:edp}). 

\subsection{Empirical Bellman Operator}\label{sec:ebo}

Recall the  definition of   the empirical Bellman operator in equation \eqref{eq:empiricalBellman-defn}. Here we give a  mathematical basis for that definition which will help us to analyze the convergence behaviour of EVI (since EVI can be framed  as an iteration of this operator).

The empirical Bellman operator is a random operator,
because it maps random samples to operators. Recall from Section \ref{sec:probfp} that we define the random operator on  the sample space $\Omega=\left[0,1\right]^{\infty}$ where primitive uncertainties
on $\Omega$ are infinite sequences of uniform noise $\omega=\left(\xi_{i}\right)_{i\geq0}$
where each $\xi_{i}$ is an independent uniform random variable on
$\left[0,1\right]$. This convention, rather than just defining $\Omega=\left[0,1\right]^{n}$
for a fixed $n\geq1$, makes convergence statements with respect to
$n$ easier to make.

%

Classical value iteration is performed by iterating the Bellman operator
$T$. Our EVI algorithm is performed by choosing $n$ and then iterating
the random operator $\widehat{T}_{n}$. So we follow the notations introduced in Section \ref{sec:probfp} and the $k$th iterate of EVI, $\hat{v}^{k}_{n}$ is given by $\hat{v}^{k}_{n} = \widehat{T}^{k}_{n} \hat{v}^{0}_{n}$ where  $\hat{v}_{n}^{0}\in\mathbb{R}^{|\mathbb{S}|}$ be an initial seed
for EVI. 

We first show that the empirical Bellman operator satisfies the Assumptions \ref{ass:fixedpoint} - \ref{ass:pn}. Then the analysis follows the results of Section \ref{sec:randomops}. 

\begin{proposition}
\label{prop:empB-ass}
The Bellman operator $T$ and the empirical Bellman operator $\widehat{T}_{n}$ (defined in equation \eqref{eq:empiricalBellman-defn}) satisfy  Assumptions \ref{ass:fixedpoint} - \ref{ass:pn}
\end{proposition}
Proof is given in Appendix \ref{pf:prop:empB-ass}.


We note that we can explicitly give an expression for $p_{n}(\epsilon)$ (of Assumption \ref{ass:pn}) as below. For proof, refer to Proposition  \ref{prop:empB-ass}: 
\begin{equation}
\label{eq:defn-pn}
 P\left\{ \|\widehat{T}_{n}v-T\, v\| < \epsilon\right\}  > p_{n}(\epsilon) := 1 - 2\,|\mathbb{K}|\, e^{-2\,(\epsilon/\alpha)^{2}n/\left(2\,\kappa^{*}\right)^{2}}.
\end{equation}
Also we note that we can also give an explicit expression for $\kappa^{*}$ of in Assumption \ref{ass:bddness} as
\begin{equation}
\label{eq:defn-kappa*}
\kappa^{*}\triangleq\frac{\max_{\left(s,a\right)\in\mathbb{K}}|c\left(s,a\right)|}{1-\alpha}.
\end{equation}
For proof, refer to the proof of Proposition  \ref{prop:empB-ass}.

\subsection{Empirical Value Iteration}\label{sec:evi-results}

Here we use the results of Section \ref{sec:randomops} for analyzing the convergence of EVI. We first give an asymptotic result. 

\begin{proposition}
\label{prop:samplecomplexity} 
For any $\delta_{1}\in\left(0,1\right)$ select $n$ such that
\[n \geq \frac{2\left(\kappa^{*}\right)^{2}}{(\epsilon_{g}/\alpha)^{2}} \log \frac{ 2|\mathbb{K}|}{\delta_{1}}\]
then, \[\limsup_{k\rightarrow\infty}\mathcal{P}\left\{ \|\hat{v}_{n}^{k}-v^{*}\| > \eta^{*}\epsilon_{g}\right\} \leq 1 - \mu_{n}(\eta^{*}) \leq \delta_{1}\]
\end{proposition}
\proof{Proof:}
From Corollary \ref{cor:probabilistic}, $\limsup_{k\rightarrow\infty}\mathcal{P}\left\{ \|\hat{v}_{n}^{k}-v^{*}\|> \eta^{*}\epsilon_{g}\right\}  \leq\mathcal{Q}\left\{ Y_{n} > \eta^{*}\right\} = 1 - \mu_{n}\left(\eta^{*}\right)$.
For $1 - \mu_{n}\left(\eta^{*}\right)$ to be less that $\delta_{1}$, we compute $n$ using Lemma \ref{lem:steadystate} as,
\[1-\delta_{1} \leq \mu_{n}(\eta^{*}) = p_{n}^{N^{*}-\eta^{*}-1} \leq p_{n} =1-2\,|\mathbb{K}|\, e^{-2\,(\epsilon_{g}/\alpha)^{2}n/\left(2\,\kappa^{*}\right)^{2}}.  \]
Thus, we get the desired result. 
\hfill \Halmos \endproof

We cannot iterate $\widehat{T}_{n}$ forever so we need a guideline
for a finite choice of $k$. This question can be answered in terms
of mixing times. The total variation distance between two probability
measures $\mu$ and $\nu$ on $\mathbb{S}$ is
\[
\|\mu-\nu\|_{TV}=\max_{S\subset\mathbb{S}}|\mu\left(S\right)-\nu\left(S\right)|=\frac{1}{2}\sum_{s\in\mathbb{S}}|\mu\left(s\right)-\nu\left(s\right)|.
\]
Let $Q_{n}^{k}$  be the marginal distribution of $Y_{n}^{k}$
on $\mathbb{N}$ at stage $k$ and
\[
d\left(k\right)=\|Q_{n}^{k}-\mu_{n}\|_{TV}
\]
be the total variation distance between $Q_{n}^{k}$ and the steady
state distribution $\mu_{n}$. For $\delta_{2} > 0$, we define
\[
t_{mix}\left(\delta_{2}\right)=\min\left\{ k\mbox{ : }d\left(k\right)\leq\delta_{2}\right\} 
\]
to be the minimum length of time needed for the marginal distribution
of $Y_{n}^{k}$ to be within $\delta_{2}$ of the steady state distribution in total variation norm. 

We now bound $t_{mix}\left(\delta_{2}\right)$.
\begin{lemma}
\label{lem:convergencerate} For any $\delta_{2}>0$, 
\[
t_{mix}\left(\delta_{2}\right)\leq\log\left(\frac{1}{\epsilon\,\mu_{n,\,\min}}\right).
\]
where $\mu_{n,\,\min} :=\min_{\eta}\mu_{n}\left(\eta\right)$. 
\end{lemma}
\proof{Proof:}
Let $\mathfrak{Q}_{n}$ be the transition matrix of the Markov chain $\{Y_{n}^{k}\}_{k \geq 0}$. Also let $\lambda_{\star}=\max\left\{ |\lambda|\mbox{ : }\lambda\mbox{ is an eigenvalue of }\mathfrak{Q}_{n},\,\lambda\ne1\right\}$. By \cite[Theorem 12.3]{Levin_Mixing_2008}, 
\[t_{mix}\left(\delta_{2}\right) \leq\log\left(\frac{1}{\delta_{2}\,\mu_{n,\,\min}}\right)\frac{1}{1-\lambda_{\star}}  \]
but $\lambda_{\star}=0$ by Lemma  given in Appendix \ref{pf:lem:eigen}.
\hfill \Halmos \endproof

We now use the above bound on mixing time to get a non-asymptotic bound for EVI.  

\begin{proposition}
\label{prop:mixingtime}
For any fixed $n \geq 1$, $\mathcal{P}\left\{ \|\hat{v}_{n}^{k}-v^{*}\| > \eta^{*}\epsilon_{g}\right\} \leq 2 \delta_{2} + (1-\mu_{n}(\eta^{*}))$ for $k \geq\log\left(\frac{1}{\delta_{2} \mu_{n,\,\min}}\right)$.
\end{proposition}
\proof{Proof:}
For $k \geq  \log\left(\frac{1}{\delta_{2} \mu_{n,\,\min}}\right) \geq  t_{\text{mix}}(\delta_{2})$, 
\[d(k) = \frac{1}{2} \sum^{N^{*}}_{i = \eta^{*}} |Q(Y^{k}_{n}=i) - \mu_{n}(i) | \leq \delta_{2}.\] Then, $ |Q(Y^{k}_{n}=\eta^{*}) - \mu_{n}(\eta^{*})| \leq 2 d(t) \leq 2\delta_{2}$. So, $\mathcal{P}\left\{ \|\hat{v}_{n}^{k}-v^{*}\| > \eta^{*}\epsilon_{g}\right\}  \leq Q(Y^{k}_{n} > \eta^{*}) = 1 - Q(Y^{k}_{n} = \eta^{*}) \leq 2 \delta_{2} + (1 - \mu_{n}(\eta^{*}))$. 
\hfill \Halmos \endproof

We now combine Proposition \ref{prop:samplecomplexity} and \ref{prop:mixingtime} to prove  Theorem \ref{thm:EVI_main}.\\ \\
\noindent{\bf Proof of Theorem \ref{thm:EVI_main}:}\\
\proof{Proof:}
Let $\epsilon_{g}=\epsilon/\eta^{*}$, and  $\delta_{1}, \delta_{2}$ be positive  with $\delta_{1} + 2 \delta_{2} \leq \delta$.  By Proposition \ref{prop:samplecomplexity}, for $n \geq n(\epsilon, \delta)$ we have
 \[\limsup_{k\rightarrow\infty}\mathcal{P}\left\{ \|\hat{v}_{n}^{k}-v^{*}\| > \epsilon \right\}  =  \limsup_{k\rightarrow\infty}\mathcal{P}\left\{ \|\hat{v}_{n}^{k}-v^{*}\| > \eta^{*}\epsilon_{g}\right\} = \leq 1 - \mu_{n}(\eta^{*}) \leq \delta_{1}.\]

Now, for $k \geq k(\epsilon, \delta)$, by Proposition \ref{prop:mixingtime},  $\mathcal{P}\left\{ \|\hat{v}_{n}^{k}-v^{*}\| > \eta^{*}\epsilon_{g}\right\} = \mathcal{P}\left\{ \|\hat{v}_{n}^{k}-v^{*}\| > \epsilon\right\} \leq 2 \delta_{2} + (1-\mu_{n}(\eta^{*}))$. Combining both we get, $\mathcal{P}\left\{ \|\hat{v}_{n}^{k}-v^{*}\| > \epsilon\right\} \leq \delta$.
\hfill \Halmos \endproof

\subsection{Empirical Policy Iteration}\label{sec:epi-results}

We now consider empirical policy iteration. EPI is different
from EVI, and seemingly more difficult to analyze, because it does
not correspond to iteration of a random operator. Furthermore, it has  two simulation components, empirical policy evaluation and empirical policy update. However, we show that the convergence analysis in a manner  similar to that of EVI. 

 We first give a sample complexity result for policy evaluation. For a policy $\pi$, let $v^{\pi}$ be the actual value of the policy and let  $\hat{v}_{q}^{\pi}$ be the empirical evaluation.  Then,

\begin{proposition}
\label{prop:simulation} 
For any $\pi\in\Pi$, $\epsilon \in (0, \gamma)$ and for any $\delta > 0$
\[
P\left\{ \|\hat{v}_{q}^{\pi}-v^{\pi}\|\geq\epsilon\right\} \leq \delta,
\]
for 
\[q \geq \frac{2(\kappa^{*} (\mathfrak{T} +1))^{2}}{ (\epsilon-\gamma)^{2}} \log \frac{2 |\mathbb{S}|}{\delta},  \]
where $\hat{v}_{q}$ is evaluation of $v^{\pi}$ by averaging $q$
simulation runs.
\end{proposition}
\proof{Proof:}
Let $v^{\pi, \mathfrak{T}} := \mathbb{E}\left[ \sum_{t=0}^{\mathfrak{T}}\alpha^{t}c\left(s_{t}\left(\boldsymbol{\omega}\right)),\pi\left(s_{t}\left(\boldsymbol{\omega}\right)\right)\right) \right]$. Then,
\begin{align*}
|\hat{v}_{q}^{\pi}(s)-v^{\pi}(s)| &\leq |\hat{v}_{q}^{\pi}(s) - v^{\pi, \mathfrak{T}}| + |v^{\pi, \mathfrak{T}} - v^{\pi}|\\
&\leq |\frac{1}{q} \sum^{q}_{i=1} \sum^{\mathfrak{T}}_{t=0} \alpha^{k}c\left(s_{t}\left(\boldsymbol{\omega}_{i}\right),\pi\left(s_{t}\left(\boldsymbol{\omega}_{i}\right)\right)\right) - v^{\pi, \mathfrak{T}}| + \gamma \\
&\leq \sum^{\mathfrak{T}}_{t=0} \left|\frac{1}{q}\sum^{q}_{i=1}\left(c\left(s_{t}\left(\boldsymbol{\omega}_{i}\right),\pi\left(s_{t}\left(\boldsymbol{\omega}_{i}\right)\right)\right) -  \mathbb{E}\left[c\left(s_{t}\left(\boldsymbol{\omega}\right),\pi\left(s_{t}\left(\boldsymbol{\omega}\right)\right)\right) \right] \right) \right| + \gamma.
\end{align*}
Then, with $\tilde{\epsilon} = (\epsilon-\gamma)/(\mathfrak{T} +1)$,

\begin{align*}
P\left(|\hat{v}_{q}^{\pi}(s)-v^{\pi}(s)| \geq \epsilon \right) &\leq  P\left(\left|\frac{1}{q}\sum^{q}_{i=1}\left(c\left(s_{t}\left(\boldsymbol{\omega}_{i}\right),\pi\left(s_{t}\left(\boldsymbol{\omega}_{i}\right)\right)\right) -  \mathbb{E}\left[c\left(s_{t}\left(\boldsymbol{\omega}_{i}\right),\pi\left(s_{t}\left(\boldsymbol{\omega}_{i}\right)\right)\right) \right] \right) \right| \geq \tilde{\epsilon} \right) \\
& \leq 2 e^{-2q\tilde{\epsilon}^{2}/(2\kappa^{*})^{2} }.
\end{align*}
By applying the union bound, we get
\[P\left(\|\hat{v}_{q}^{\pi}-v^{\pi}\|\geq\epsilon\right) \leq  2 |\mathbb{S}| e^{-2q\tilde{\epsilon}^{2}/(2\kappa^{*})^{2} }  .\]
For $q \geq \frac{2(\kappa^{*} (\mathfrak{T} +1))^{2}}{ (\epsilon-\gamma)^{2}} \log \frac{2 |\mathbb{S}|}{\delta} $ the above probability is less than $\delta$. 
\hfill \Halmos \endproof

We define 
\begin{equation}
\label{eq:defn-rq}
P\left(\|\hat{v}_{q}^{\pi}-v^{\pi}\| < \epsilon\right) > r_{q}(\epsilon) : = 1 - 2 |\mathbb{S}| e^{-2q\tilde{\epsilon}^{2}/(2\kappa^{*})^{2} },~\text{with}~\tilde{\epsilon} = (\epsilon-\gamma)/(\mathfrak{T} +1).
\end{equation}

We say that  empirical policy evaluation is $\epsilon_{1}$-accurate if $\|\hat{v}_{q}^{\pi}-v^{\pi}\| < \epsilon_{1}$. Then by the above proposition empirical policy evaluation is $\epsilon_{1}$-accurate with a probability greater than $r_{q}(\epsilon_{1})$. 

The accuracy of empirical policy update compared to the actual policy update indeed depends on the empirical Bellman operator $\widehat{T}_{n}$. We say that empirical policy update is $\epsilon_{2}$-accurate if $\|\widehat{T}_{n}v - Tv\| < \epsilon_{2}$. Then, by the definition of  $p_{n}$ in equation \eqref{eq:defn-pn}, our empirical policy update is $\epsilon_{2}$-accurate with a probability greater than $p_{n}(\epsilon_{2})$ 

We now give an important technical lemma. Proof is essentially a probabilistic modification of Lemmas 6.1 and 6.2 in \cite{Bertsekas_Neuro_1996} and is omitted.

\begin{lemma}
\label{lem:EPI_error} 
Let $\{\pi_k\}_{k \geq 0}$ be the sequence of policies from the EPI algorithm. For a fixed $k$, assume that $P\left(\|v^{\pi_{k}} - \hat{v}_{q}^{\pi_{k}}\| < \epsilon_{1} \right) \geq (1 - \delta_{1})$ and $P\left(\|T\hat{v}_{q}^{\pi_{k}} - \widehat{T}_{n}\hat{v}_{q}^{\pi_{k}}\| < \epsilon_{2}\right) \geq (1 - \delta_{2})$.  Then,
\[\|v^{\pi_{k+1}} - v^{*} \| \leq \alpha \|v^{\pi_{k}} - v^{*} \| + \frac{\epsilon_{2}+2\alpha \epsilon_{1}}{(1-\alpha)}  \]
with probability at least  $(1 - \delta_{1}) (1 - \delta_{2}).$
\end{lemma}

We now proceed as in the analysis of EVI given in the previous subsection. Here we track the sequence  $\{\|v^{\pi_{k}} - v^{*} \|\}_{k \geq 0}$. Note that this being a proof technique, the fact that  the value $\|v^{\pi_{k}} - v^{*} \|$ is not observable does not affect our algorithm or its convergence behavior. We define   
\[X_{n,\, q}^{k} =\left\lceil \|\hat{v}^{\pi_{k}}-v^{*}\|/\epsilon_{g} \right\rceil   \]
 where the granularity $\epsilon_{g}$ is fixed according to the problem parameters as $\epsilon_{g} = \frac{\epsilon_{2}+2\alpha \epsilon_{1}}{(1-\alpha)}$. Then by Lemma \ref{lem:EPI_error}, 
 \begin{equation}
\label{eq:xk-to-yk-step-epi}
\text{if}~X_{n,\, q}^{k}=\eta,~\text{then}~X_{n,\, q}^{k+1} \leq \left\lceil \alpha\,\eta+1\right\rceil ~\text{with a probability at least}~ p_{n,q}=r_{q}(\epsilon_{1}) p_{n}(\epsilon_{2}).\
\end{equation}

This is equivalent to the result for EVI given in equation \eqref{eq:xk-to-yk-step-evi}. Hence the analysis is the same from here onwards.  However, for completeness, we explicitly give the dominating Markov chain and its steady state distribution.

For $p_{n,q}$ given in display \eqref{eq:xk-to-yk-step-epi}, we construct the  dominating Markov chain  $\left\{ Y_{n,\, q}^{k}\right\} _{k\geq0}$ as
\begin{equation}
Y_{n,\, q}^{k+1}=\begin{cases}
\max\left\{ Y_{n,\, q}^{k}-1,\,\eta^{*}\right\} , & \mbox{w.p. }p_{n,\, q},\\
N^{*}, & \mbox{w.p. }1-p_{n,\, q},
\end{cases}
\end{equation}
which exists on the state space $\mathcal{X}$. 
The family $\left\{ Y_{n,\, q}^{k}\right\} _{k\geq0}$ is identical
to $\left\{ Y_{n}^{k}\right\} _{k\geq0}$ except that its transition
probabilities depend on $n$ and $q$ rather than just $n$.  Let $\mu_{n,q}$ denote the steady state distribution of the Markov chain $\left\{ Y_{n,q}^{k}\right\} _{k\geq0}$. Then by  Lemma \ref{lem:steadystate},
\begin{equation}
\label{eq:mu-nq}
\mu_{n,q}\left(\eta^{*}\right)=\,~p_{n,q}^{N^{*}-\eta^{*}-1}, \mu_{n,q}\left(N^{*}\right)=\,  \frac{1-p_{n,q}}{p_{n,q}},~ \mu_{n,q}\left(i\right)=\,  \left(1-p_{n,q}\right)p_{n,q}^{\left(N^{*}-i-1\right)},  \forall i=\eta^{*}+1,\ldots,N^{*}-1 . 
\end{equation}
\\\noindent{\bf Proof of Theorem \ref{thm:EPI_main}:}
\proof{Proof:}
First observe that by the given choice of $n$ and $q$, we have $r_{q} \geq (1-\delta_{11})$ and $p_{n} \geq (1-\delta_{12})$. Hence $1 - p_{n,q} \leq \delta_{11} + \delta_{12} - \delta_{11}\delta_{12} < \delta_{1}$. 

Now by Corollary \ref{cor:probabilistic},
\[ \limsup_{k\rightarrow\infty}\mathcal{P}\left\{ \|\hat{v}^{\pi_{k}}-v^{*}\| > \epsilon \right\} = \limsup_{k\rightarrow\infty}\mathcal{P}\left\{ \|\hat{v}^{\pi_{k}}-v^{*}\| > \eta^{*}\epsilon_{g}\right\}   \leq\mathcal{Q}\left\{ Y_{n,q} > \eta^{*}\right\} = 1 - \mu_{n,q}\left(\eta^{*}\right).\]
For $ 1 - \mu_{n,q}\left(\eta^{*}\right)$ to be less than $\delta_{1}$ we need $1-\delta_{1} \leq \mu_{n}(\eta^{*}) = p_{n,q}^{N^{*}-\eta^{*}-1} \leq p_{n,q}$ which true as verified above. Thus we get
\[ \limsup_{k\rightarrow\infty}\mathcal{P}\left\{ \|\hat{v}^{\pi_{k}}-v^{*}\| > \epsilon \right\} \leq \delta_{1},\]
similar to the result of Proposition \ref{prop:samplecomplexity}. Selecting the number of iterations $k$ based on the mixing time is same as given in Proposition \ref{prop:mixingtime}. Combining both as in the proof of Theorem \ref{thm:EVI_main} we get the desired result. 
\hfill \Halmos \endproof

\section{Variations and Extensions}\label{sec:variations}

We now consider some variations and extensions of EVI.

\subsection{Asynchronous Value Iteration}\label{sec:asynch}

The EVI algorithm described above is synchronous, meaning that the
value estimates for every state are updated simultaneously. Here we consider each state 
to be visited at least once to complete a full update cycle. We modify
the earlier argument to account for the possibly random time between
full update cycles.

Classical asynchronous value iteration with exact updates has already
been studied. Let $\left(x_{k}\right)_{k\geq0}$ be any infinite sequence
of states in $\mathbb{S}$. This sequence $\left(x_{k}\right)_{k\geq0}$
may be deterministic or stochastic, and it may even depend online
on the value function updates. For any $x\in\mathbb{S}$, we define
the asynchronous Bellman operator $T_{x}\mbox{ : }\mathbb{R}^{|\mathbb{S}|}\rightarrow\mathbb{R}^{|\mathbb{S}|}$
via
\[
\left[T_{x}v\right]\left(s\right)=\begin{cases}
\min_{a\in A\left(s\right)}\left\{ c\left(s,a\right)+ \alpha \mathbb{E}\left[v\left(\psi\left(s,a,\xi\right)\right)\right]\right\} , & s=x,\\
v\left(s\right), & s\ne x.
\end{cases}
\]
The operator $T_{x}$ only updates the estimate of the value function
for state $x$, and leaves the estimates for all other states exactly
as they are. Given an initial seed $v^{0}\in\mathbb{R}^{|\mathbb{S}|}$,
asynchronous value iteration produces the sequence $\left\{ v^{k}\right\} _{k\geq0}$
defined by $v^{k}=T_{x_{t}}T_{x_{t-1}}\cdots T_{x_{0}}v^{0}$ for
$k\geq0$.

The following key properties of $T_{x}$ are immediate.  
\begin{lemma}
\label{lem:asyn-Tx}
For any $x\in\mathbb{S}$:

(i) $T_{x}$ is monotonic;

(ii) $T_{x}\left[v+\eta\,1\right]=T_{x}v+\alpha\,\eta\, e_{x}$, where $e_{x}\in\mathbb{R}^{|\mathbb{S}|}$
be the unit vector corresponding to $x\in\mathbb{S}$.

\end{lemma}
Proof is given in Appendix \ref{pf:lem:asyn-Tx}

The next lemma is used to show that classical asynchronous VI converges.
Essentially, a cycle of updates that visits every state at least once
is a contraction. 
\begin{lemma}
Let $\left(x_{k}\right)_{k=1}^{K}$ be any finite sequence of states
such that every state in $\mathbb{S}$ appears at least once, then
the operator
\[
\widetilde{T}=T_{x_{1}}T_{x_{2}}\cdots T_{x_{K}}
\]
is a contraction with constant $\alpha$.
\end{lemma}
It is known that asynchronous VI converges when each state is visited
infinitely often. To continue, define $K_{0}=0$ 
and in general, we define
\[
K_{m+1}\triangleq\inf\left\{ k\mbox{ : }k\geq K_{m},\,\left(x_{i}\right)_{i=K_{m}+1}^{k}\mbox{ includes every state in }\mathbb{S}\right\} .
\]
Time $K_{1}$ is the first time that every state in $\mathbb{S}$
is visited at least once by the sequence $\left(x_{k}\right)_{k\geq0}$.
Time $K_{2}$ is the first time after $K_{1}$ that every state is
visited at least once again by the sequence $\left(x_{k}\right)_{k\geq0}$,
etc. The times $\left\{ K_{m}\right\} _{m\geq0}$ completely depend
on $\left(x_{k}\right)_{k\geq0}$. For any $m\geq0$, if we define
\[
\widetilde{T}=T_{K_{m+1}}T_{K_{m+1}-1}\cdots T_{K_{m}+2}T_{K_{m}+1},
\]
then we know
\[
\|\widetilde{T}\, v-v^{*}\|\leq\alpha\,\|v-v^{*}\|,
\]
by the preceding lemma. It is known that asynchronous VI converges
under some conditions on $\left(x_{k}\right)_{k\geq0}$.
\begin{theorem}\cite{bertsekas1995dynamic}. Suppose each state in $\mathbb{S}$ is included infinitely often by
the sequence $\left(x_{k}\right)_{k\geq0}$. Then $v^{k}\rightarrow v^{*}$.
\end{theorem}
Next we describe an empirical version of classical asynchronous value
iteration. Again, we replace exact computation of the expectation
with an empirical estimate, and we regenerate the sample at each iteration.

\begin{algorithm}[!tph]
\caption{Asynchronous empirical value iteration}

Input: $v^{0}\in\mathbb{R}^{|\mathbb{S}|}$, sample size $n\geq1$,
a sequence $\left(x_{k}\right)_{k\geq0}$.

Set counter $k=0$.
\begin{enumerate}
\item Sample $n$ uniformly distributed random variables $\left\{ \xi_{i}\right\} _{i=1}^{n}$,
and compute
\[
\hat{v}^{k+1}\left(s\right)=\begin{cases}
\min_{a\in A\left(s\right)}\left\{ c\left(s,a\right)+\frac{\alpha}{n}\sum_{i=1}^{n}\hat{v}^{k}\left(\psi\left(s,a,\xi_{i}\right)\right)\right\} , & s=x_{k},\\
v\left(s\right), & s\ne x_{k}.
\end{cases}
\]

\item Increment $k:=k+1$ and return to step 2.\end{enumerate}
\end{algorithm}

Step 1 of this algorithm replaces the exact computation $v^{k+1}=T_{x_{k}}v^{k}$
with an empirical variant. Using our earlier notation, we let $\widehat{T}_{x,n}$
be a random operator that only updates the value function for state
$x$ using an empirical estimate with sample size $n\geq1$:
\[
\left[\widehat{T}_{x,n}\left(\omega\right)v\right]\left(s\right)=\begin{cases}
\min_{a\in A\left(s\right)}\left\{ c\left(s,a\right)+\frac{1}{n}\sum_{i=1}^{n}v\left(\psi\left(s,a,\xi_{i}\right)\right)\right\} , & s=x,\\
v\left(s\right), & s\ne x.
\end{cases}
\]
We use $\left\{ \hat{v}_{n}^{k}\right\} _{k\geq0}$ to denote the
sequence of asynchronous EVI iterates, 
\[
\hat{v}_{n}^{k+1}=\widehat{T}_{x_{k},n}\widehat{T}_{x_{k-1},n}\cdots\widehat{T}_{x_{0},n}\hat{v}_{n}^{0},
\]
or more compactly $\hat{v}_{n}^{k+1}=\widehat{T}_{x_{k},n}\hat{v}_{n}^{k}$
for all $k\geq0$.

We can use a slightly modified stochastic dominance argument to show
that asynchronous EVI converges in a probabilistic sense. Only now
we must account for the hitting times $\left\{ K_{m}\right\} _{m\geq0}$
as well, since the accuracy of the overall update depends on the accuracy
in $\widehat{T}_{x,n}$ as well as the length of the interval $\left\{ K_{m}+1,\, K_{m}+2,\ldots,K_{m+1}\right\} $.
In asynchronous EVI, we will focus on $\left\{ \hat{v}_{n}^{K_{m}}\right\} _{m\geq0}$
rather than $\left\{ \hat{v}_{n}^{k}\right\} _{k\geq0}$. We  check the progress of the algorithm at the end of complete
update cycles.

In the simplest update scheme, we can order the states and then update
them in the same order throughout the algorithm. The set $\left(x_{k}\right)_{k\geq0}$
is deterministic in this case, and the intervals $\left\{ K_{m}+1,\, K_{m}+2,\ldots,K_{m+1}\right\} $
all have the same length $|\mathbb{S}|$. Consider
\[
\widetilde{T}=T_{x_{K_{1}},n}T_{x_{K_{1}-1},n}\cdots T_{x_{1},n}T_{x_{0},n},
\]
the operator $\widehat{T}_{x_{0},n}$ introduces $\epsilon$ error
into component $x_{0}$,  the operator $\widehat{T}_{x_{1},n}$
introduces $\epsilon$ error into component $x_{1}$, etc. To ensure
that
\[
\widehat{T}=\widehat{T}_{x_{K_{1}},n}\widehat{T}_{x_{K_{1}-1},n}\cdots\widehat{T}_{x_{1},n}\widehat{T}_{x_{0},n}
\]
is close to $\widetilde{T}$, we require each $\widehat{T}_{x_{k},n}$
to be close to $T_{x_{k}}$ for $k =  0,\,1,\ldots,K_{1}-1,K_{1}$.

The following noise driven perspective helps with our error analysis.
In general, we can view asynchronous empirical value iteration as
\[
v'=T_{x}\, v+\varepsilon
\]
for all $k\geq0$ where
\[
\varepsilon=\widehat{T}_{x,n}v-T_{x}v
\]
is the noise for the evaluation of $T_{x}$ (and it has at most one
nonzero component).

Starting with $v^{0}$, define the sequence $\left\{ v^{k}\right\} _{k\geq0}$
by exact asynchronous value iteration $v^{k+1}=T_{x_{k}}v^{k}$ for
all $k\geq0$. Also set $\tilde{v}_{0}:=v_{0}$ and define
\[
\tilde{v}^{k+1}=T_{x_{k}}\tilde{v}^{k}+\varepsilon_{k}
\]
for all $k\geq0$ where $\varepsilon_{k}\in\mathbb{R}^{|\mathbb{S}|}$
is the noise for the evaluation of $T_{x_{k}}$ on $\tilde{v}^{k}$.
In the following proposition, we compare the sequences of value functions
$\left\{ v^{k}\right\} _{k\geq0}$ and $\left\{ \tilde{v}^{k}\right\} _{k\geq0}$
under conditions on the noise $\left\{ \varepsilon_{k}\right\} _{k\geq0}$.

\begin{proposition}
\label{prop:Asynchronous_ErrorBound} 
Suppose $-\eta\,1\leq\varepsilon_{i}\leq\eta\,1$
for all $j=0,1,\ldots,k$ where $\eta\geq0$ and $1\in\mathbb{R}^{|\mathbb{S}|}$,
i.e. the error is uniformly bounded for $j=0,1,\ldots,k$. Then, for
all $j=0,1,\ldots,k$: 
\[
v_{j}-\left(\sum_{i=0}^{j}\alpha^{i}\right)\eta\,1\leq\tilde{v}_{j}\leq v_{j}+\left(\sum_{i=0}^{j}\alpha^{i}\right)\eta\,1.
\]
\end{proposition}
Proof is given in Appendix \ref{pf:prop:Asynchronous_ErrorBound}

Now we can use the previous proposition to obtain conditions for $\|\widetilde{T}\, v-\widehat{T}\, v\|<\epsilon$
(for our deterministic update sequence). Starting with the update
for state $x_{0}$, we can choose $n$ to ensure
\[
\|T_{x_{0},n}v-\widehat{T}_{x_{0},n}v\|<\epsilon/|\mathbb{S}|
\]
similar to that in equation \eqref{eq:defn-pn}. However,
in this case our error bound is
\begin{align*}
P\left\{ \|\widehat{T}_{x,n}v-T_{x}v\|\geq\epsilon/|\mathbb{S}|\right\} \leq\, & P\left\{ \max_{a\in A\left(s\right)}|\frac{1}{n}\sum_{i=1}^{n}v\left(\psi\left(s,a,\xi_{i}\right)\right)-\mathbb{E}\left[v\left(\psi\left(s,a,\xi\right)\right)\right]|\geq\epsilon/(\alpha |\mathbb{S}|)\right\} \\
\leq\, & 2\,|\mathbb{A}|\, e^{-2\,\left(\epsilon/(\alpha |\mathbb{S}|)\right)^{2}n/\left(2\,\kappa^{*}\right)^{2}},
\end{align*}
for all $v\in\mathbb{R}^{|\mathbb{S}|}$ (which does not depend on
$x$). We are only updating one state, so we are concerned with the
approximation of at most $|\mathbb{A}|$ terms $c\left(s,a\right)+ \alpha \mathbb{E}\left[v\left(\psi\left(s,a,\xi\right)\right)\right]$
rather than $|\mathbb{K}|$. At the next update we want

\[
\|\widehat{T}{}_{x_{1},n}\hat{v}_{n}^{1}-T_{x_{1},n}\hat{v}_{n}^{1}\|<\epsilon/|\mathbb{S}|,
\]
and we get the same error bound as above.

Based on this reasoning, assume
\[
\|T_{x_{k},n}\hat{v}_{n}^{k}-\widehat{T}_{x_{k},n}\hat{v}_{n}^{k}\|<\epsilon/|\mathbb{S}|
\]
for all $k=0,1,\ldots,K_{1}$. In this case we will in fact get the
stronger error guarantee
\[
\|\widetilde{T}\, v-\widehat{T}\, v\|<\frac{\epsilon}{|\mathbb{S}|}\sum_{i=0}^{|\mathbb{S}|-1}\alpha^{i}<\epsilon
\]
from Proposition \ref{prop:Asynchronous_ErrorBound}. The complexity
estimates are multiplicative, so the probability $\|T_{x_{k},n}\hat{v}_{n}^{k}-\widehat{T}_{x_{k},n}\hat{v}_{n}^{k}\|<\epsilon/|\mathbb{S}|$
for all $k=0,1,\ldots,K_{1}$ is bounded above by
\[
p_{n}=2\,|\mathbb{S}|\,|\mathbb{A}|\, e^{-2\,\left(\epsilon/(\alpha |\mathbb{S}|)\right)^{2}n/\left(2\,\kappa^{*}\right)^{2}}.
\]
To understand this result, remember that $|\mathbb{S}|$ iterations
of asynchronous EVI amount to at most $|\mathbb{S}|\,|\mathbb{A}|$
empirical estimates of $c\left(s,a\right)+\mathbb{E}\left[v\left(\psi\left(s,a,\xi\right)\right)\right]$.
We require all of these estimates to be within error $\epsilon/|\mathbb{S}|$.

We can take the above value for $p_{n}$ and apply our earlier stochastic
dominance argument to\\
$\left\{ \|\hat{v}_{n}^{K_{m}}-v^{*}\|\right\} _{m\geq0}$, without
further modification. This technique extends to any deterministic
sequence $\left(x_{k}\right)_{k\geq0}$ where the lengths of a full
update for all states $|K_{m+1}-K_{m}|$ are uniformly bounded for
all $m\geq0$ (with the sample complexity estimate suitably adjusted).

\subsection{Minimax Value Iteration}\label{sec:minimax}

Now we consider a two player zero sum Markov game and show how an
empirical min-max value iteration algorithm can be used to a compute
an approximate Markov Perfect equilibrium. Let the Markov game be described by
the 7-tuple
\[
\left(\mathbb{S},\,\mathbb{A},\,\left\{ A\left(s\right):s\in\mathbb{S}\right\} ,\,\mathbb{B},\,\left\{ B\left(s\right):s\in\mathbb{S}\right\} ,\, Q,\, c\right).
\]
The action space $\mathbb{B}$ for player 2 is finite and $B\left(s\right)$
accounts for feasible actions for player 2. We let
\[
\mathbb{K}=\left\{ \left(s,a,b\right):s\in\mathbb{S},\, a\in A\left(s\right),\, b\in B\left(s\right)\right\} 
\]
be the set of feasible station-action pairs. The transition law $Q$
governs the system evolution, $Q\left(\cdot|s,a,b\right)\in\mathcal{P}\left(\mathbb{A}\right)$
for all $\left(s,a\right)\in\mathbb{K}$, which is the probability of next visiting the state
$j$ given $\left(s,a,b\right)$. Finally, $c:\mathbb{K}\rightarrow\mathbb{R}$
is the cost function (say of player 1) in state $s$ for actions $a$ and
$b$. Player 1 wants to minimize this quantity, and player 2 is
trying to maximize this quantity.

Let the operator $T$ be defined as  $T:\mathbb{R}^{|\mathbb{S}|}\rightarrow\mathbb{R}^{|\mathbb{S}|}$
is defined as
\[
\left[T\, v\right]\left(s\right)\triangleq\min_{a\in A\left(s\right)}\max_{b\in B\left(s\right)}\left\{ c\left(s,a,b\right)+\alpha\,\mathbb{E}\left[v\left(\tilde{s}\right)|s,a,b\right]\right\} ,\,\forall s\in\mathbb{S},
\]
for any $v\in\mathbb{R}^{|\mathbb{S}|}$, where $\tilde{s}$ is the
random next state visited and
\[
\mathbb{E}\left[v\left(\tilde{s}\right)|s,a,b\right]=\sum_{j\in\mathbb{S}}v\left(j\right)Q\left(j|s,a,b\right)
\]
is the same  expected cost-to-go for player
1. We call $T$ the Shapley operator in honor of Shapley who first introduced it \cite{shapley1953stochastic}. We can use $T$ to compute the optimal value function of the same which is given by
\[
v^{*}\left(s\right)=\max_{a\in A\left(s\right)}\min_{b\in B\left(s\right)}\left\{ c\left(s,a,b\right)+\alpha\,\mathbb{E}\left[v^{*}\left(\tilde{s}\right)|s,a,b\right]\right\} ,\,\forall s\in\mathbb{S},
\]
is the optimal value function for player 1.

It is well known that that the  Shapley operator is a contraction mapping.
\begin{lemma}
\label{lem:ShapleyT}
The Shapley operator $T$ is a contraction.
\end{lemma}
Proof is given in Appendix \ref{pf:lem:ShapleyT} for completeness.

To compute $v^{*}$, we can iterate $T$. Pick any initial seed $v^{0}\in\mathbb{R}^{\mathbb{S}}$,
take $v^{1}=T\, v^{0}$, $v^{2}=T\, v^{1}$, and in general $v^{k+1}=T\, v^{k}$
for all $k\geq0$. It is known that \cite{shapley1953stochastic} this procedure converges to the
optimal value function. We refer to this  as the classical minimax  value iteration.

Now, using the simulation model $\psi:\mathbb{S}\times\mathbb{A}\times\mathbb{B}\times\left[0,1\right]\rightarrow\mathbb{S}$,
the empirical Shapley  operator can be written as
\[
\left[T\, v\right]\left(s\right)\triangleq\max_{a\in A\left(s\right)}\min_{b\in B\left(s\right)}\left\{ c\left(s,a,b\right)+\alpha\,\mathbb{E}\left[v\left(\psi\left(s,a,b,\xi\right)\right)\right]\right\} ,\,\forall s\in\mathbb{S},
\]
where $\xi$ is a uniform random variable on $\left[0,1\right]$.

We will replace the expectation $\mathbb{E}\left[v\left(\psi\left(s,a,b,\xi\right)\right)\right]$
with an empirical estimate. Given a sample of $n$ uniform random
variables, $\left\{ \xi_{i}\right\} _{i=1}^{n}$, the empirical estimate
of $\mathbb{E}\left[v\left(\psi\left(s,a,b,\xi\right)\right)\right]$
is $\frac{1}{n}\sum_{i=1}^{n}v\left(\psi\left(s,a,b,\xi_{i}\right)\right)$.
Our algorithm is summarized next.

\begin{algorithm}[!tph]
\caption{Empirical value iteration for minimax}

Input: $v^{0}\in\mathbb{R}^{\mathbb{S}}$, sample size $n\geq1$.

Set counter $k=0$.
\begin{enumerate}
\item Sample $n$ uniformly distributed random variables $\left\{ \xi_{i}\right\} _{i=1}^{n}$,
and compute
\[
v^{k+1}\left(s\right)=\max_{a\in A\left(s\right)}\min_{b\in B\left(s\right)}\left\{ c\left(s,a,b\right)+\frac{\alpha}{n}\sum_{i=1}^{n}v^{k}\left(\psi\left(s,a,b,\xi_{i}\right)\right)\right\} ,\,\forall s\in\mathbb{S}.
\]

\item Increment $k:=k+1$ and return to step 2.\end{enumerate}
\end{algorithm}

In each iteration, we take a new set of samples and use this empirical
estimate to approximate $T$. Since $T$ is a contraction with a known
convergence rate $\alpha$, we can apply the exact same development
as for empirical value iteration.

\subsection{The Newsvendor Problem}\label{sec:newsvendor}

We now show via the newsvendor problem that the empirical dynamic programming method can sometimes work remarkably well even for continuous states and action spaces. This, of course, exploits the linear structure of the newsvendor problem. 
%

Let $D$ be a continuous random variable representing the stationary demand distribution.
Let $\left\{ D_{k}\right\} _{k\geq0}$ be independent and identically
distributed collection of random variables with the same distribution
as $D$, where $D_{k}$ is the demand in period $k$. The unit order
cost is $c$, unit holding cost is $h$, and unit backorder cost is
$b$. We let $x_{k}$ be the inventory level at the beginning of period
$k$, and we let $q_{k}\geq0$ be the order quantity before demand
is realized in period $k$.

For technical convenience, we only allow stock levels in the compact
set $\mathcal{X}=\left[x_{\min},\, x_{\max}\right]\subset\mathbb{R}$.
This assumption is not too restrictive, since a firm would not want a large
number of backorders and any real warehouse has finite capacity. Notice
that since we restrict to $\mathcal{X}$, we know that no order quantity
will ever exceed $q_{\max}=x_{\max}-x_{\min}$. Define the continuous
function $\psi:\mathbb{R}\rightarrow\mathcal{X}$ via

\[
\psi\left(x\right)=\begin{cases}
x_{\max}, & \mbox{if }x>x_{\max},\\
x_{\min}, & \mbox{if }x<x_{\min},\\
x, & \mbox{otherwise},
\end{cases}
\]
The function $\psi$ accounts for the state space truncation. The
system dynamic is then
\begin{align*}
x_{k+1}=\psi\left(x_{k}+q_{k}-D_{k}\right),\, & \forall k\geq0.
\end{align*}
We want to solve
\begin{equation}
\inf_{\pi\in\Pi}\mathbb{E}_{\nu}^{\pi}\left[\sum_{k=0}^{\infty}\alpha^{k}\left(c\, q_{k}+\max\left\{ h\, x_{k},-b\, x_{k}\right\} \right)\right],\label{NEWSVENDOR}
\end{equation}
subject to the preceding system dynamic. We know that there is an optimal
stationary policy for this problem which only depends on the current
inventory level. The optimal cost-to-go function for this problem,
$v^{*}$, satisfies
\begin{align*}
v^{*}\left(x\right)=\inf_{q\geq0}\left\{ c\, q+\max\left\{ h\, x,-b\, x\right\} +\mathbb{E}\left[v^{*}\left(\psi\left(x+q-D\right)\right)\right]\right\} ,\, & \forall x\in\mathbb{R},
\end{align*}
where, the optimal value function $v^{*}:\mathbb{R}\rightarrow\mathbb{R}$. We will compute $v^{*}$ by
iterating an appropriate Bellman operator.

Classical value iteration for Problem (\ref{NEWSVENDOR}) consists
of iteration of an operator in $\mathcal{C}\left(\mathcal{X}\right)$,
the space of continuous functions $f:\mathcal{X}\rightarrow\mathbb{R}$.
We equip $\mathcal{C}\left(\mathcal{X}\right)$ with the norm
\[
\|f\|_{\mathcal{C}\left(\mathcal{X}\right)}=\sup_{x\in\mathcal{X}}|f\left(x\right)|.
\]
Under this norm, $\mathcal{C}\left(\mathcal{X}\right)$ is a Banach
space.

Now, the Bellman operator $T:\mathcal{C}\left(\mathcal{X}\right)\rightarrow\mathcal{C}\left(\mathcal{X}\right)$
for the newsvendor problem is given by
\[
\left[T\, v\right]\left(x\right)=\inf_{q\geq0}\left\{ c\, q+\max\left\{ h\, x,-b\, x\right\} +\alpha\,\mathbb{E}\left[v\left(\psi\left(x+q-D\right)\right)\right]\right\} ,\,\forall x\in\mathcal{X}.
\]
Value iteration for the newsvendor can then be written succinctly
as $v^{k+1}=T\, v^{k}$ for all $k\geq0$. We confirm that $T$ is
a contraction with respect to $\|\cdot\|_{\mathcal{C}\left(\mathcal{X}\right)}$
in the next lemma, and thus the Banach fixed point theorem applies.
\begin{lemma}
(i) $T$ is a contraction on $\mathcal{C}\left(\mathcal{X}\right)$
with constant $\alpha$.

(ii) Let $\left\{ v^{k}\right\} _{k\geq0}$ be the sequence produced
by value iteration, then $\lim_{k\rightarrow0}\|v^{k}-v^{*}\|_{\mathcal{C}\left(\mathcal{X}\right)}\rightarrow0$.\end{lemma}
\proof{Proof:}
(i) Choose $u,\, v\in\mathcal{C}\left(\mathcal{X}\right)$, and use
Fact \ref{fact} to compute
\begin{align*}
\|T\, u-T\, v\|_{\mathcal{C}\left(\mathcal{X}\right)}=\, & \sup_{x\in\mathcal{X}}|\left[T\, u\right]\left(x\right)-\left[T\, v\right]\left(x\right)|\\
\leq\, & \sup_{x\in\mathcal{X},\, q\in\left[0,q_{\max}\right]}\alpha\,|\mathbb{E}\left[u\left(\psi\left(x+q-D\right)\right)\right]-\mathbb{E}\left[v\left(\psi\left(x+q-D\right)\right)\right]|\\
\leq\, & \alpha\sup_{x\in\mathcal{X},\, q\in\left[0,q_{\max}\right]}\mathbb{E}\left[|u\left(\psi\left(x+q-D\right)\right)-v\left(\psi\left(x+q-D\right)\right)|\right]\\
\leq\, & \alpha\,\|u-v\|_{\mathcal{C}\left(\mathcal{X}\right)}.
\end{align*}
(ii) Since $\mathcal{C}\left(\mathcal{X}\right)$ is a Banach space
and $T$ is a contraction by part (i), the Banach fixed point theorem
applies.
\hfill \Halmos \endproof
Choose the initial form for the optimal value function as
\[
v^{0}\left(x\right)=\max\left\{ h\, x,-b\, x\right\} ,\,\forall x\in\mathcal{X}.
\]
It is chosen to represent the terminal cost in state $x$ when there
are no further ordering decisions. Then, value iteration yields 
\[
v^{k+1}\left(x\right)=\inf_{q\geq0}\left\{ c\, q+\max\left\{ h\, x,-b\, x\right\} +\alpha\,\mathbb{E}\left[v^{k}\left(\psi\left(x+q-D\right)\right)\right]\right\} ,\,\forall x\in\mathcal{X}.
\]
We note some key properties of these value functions.

\begin{lemma}
\label{lem:Lipschitz} Let $\left\{ v^{k}\right\} _{k\geq0}$ be the
sequence produced by value iteration, then $v^{k}$ is Lipschitz continuous
with constant $\max\left\{ |h|,\,|b|\right\} \sum_{i=0}^{k}\alpha^{i}$
for all $k\geq0$.\end{lemma}
\proof{Proof:}
First observe that $v^{0}$ is Lipschitz continuous with constant
$\max\left\{ |h|,\,|b|\right\} $. For $v^{1}$, we choose $x$ and
$x'$ and compute
\begin{align*}
|v^{1}\left(x\right)-v^{1}\left(x'\right)|\leq\, & \sup_{q\geq0}\left\{ \max\left\{ h\, x,-b\, x\right\} -\max\left\{ h\, x',-b\, x'\right\} \right.\\
 & \left.+\alpha\left(\mathbb{E}\left[v^{0}\left(\psi\left(x+q-D\right)\right)\right]-\mathbb{E}\left[v^{0}\left(\psi\left(x'+q-D\right)\right)\right]\right)\right\} \\
\leq\, & \max\left\{ |h|,\,|b|\right\} |x-x'|+\alpha\,\max\left\{ |h|,\,|b|\right\} \mathbb{E}\left[|\psi\left(x+q-D\right)-\psi\left(x'+q-D\right)|\right]\\
\leq\, & \max\left\{ |h|,\,|b|\right\} |x-x'|+\alpha\,\max\left\{ |h|,\,|b|\right\} |x-x'|,
\end{align*}
where we use the fact that the Lipschitz constant of $\psi$ is one.
The inductive step is similar.
\hfill \Halmos \endproof
From Lemma \ref{lem:Lipschitz}, we also conclude that the Lipschitz
constant of any iterate $v^{k}$ is bounded above by
\[
L^{*}\triangleq\max\left\{ |h|,\,|b|\right\} \sum_{i=0}^{\infty}\alpha^{i}=\frac{\max\left\{ |h|,\,|b|\right\} }{1-\alpha}.
\]
We acknowledge the dependence of Lemma \ref{lem:Lipschitz} on the
specific choice of the initial seed, $v^{0}$.

We can do empirical value iteration with the same initial seed $\hat{v}_{n}^{0}=v^{0}$ as above. Now, for $k \geq 0$,
\[
\hat{v}_{n}^{k+1}\left(x\right)=\inf_{q\geq0}\left\{ c\, q+\max\left\{ h\, x,-b\, x\right\} +\frac{\alpha}{n}\sum_{i=1}^{n}\hat{v}_{n}^{k}\left(\psi\left(x+q-D_{i}\right)\right)\right\} ,\,\forall x\in\mathbb{R}.
\]
Note that $\left\{ D_{1},\ldots,D_{n}\right\} $ is an i.i.d. sample
from the demand distribution. It is possible to perform these value
function updates exactly for finite $k$ based on \cite{Cooper_Empirical_2011}.
Also note that, the initial seed is piecewise linear with finitely many
breakpoints. Because the demand sample is finite
in each iteration, thus each iteration will take a piecewise linear function
as input and then produce a piecewise linear function as output (both
with finitely many breakpoints). Lemma \ref{lem:Lipschitz} applies
without modification to $\left\{ \hat{v}_{n}^{k}\right\} _{k\geq0}$,
all of these functions are Lipschitz continuous with constants bounded
above by $L^{*}$.

As earlier, we define the empirical Bellman operator $\widehat{T}_{n}:\Omega\rightarrow\mathcal{C}$
as
\[
\left[\widehat{T}_{n}\left(\omega\right)v\right]\left(x\right)=\inf_{q\geq0}\left\{ c\, q+\max\left\{ h\, x,-b\, x\right\} +\frac{\alpha}{n}\sum_{i=1}^{n}v\left(\psi\left(x+q-D_{i}\right)\right)\right\} ,\,\forall x\in\mathbb{R}.
\]
With the empirical Bellman operator, we write the iterates of EVI
as $\hat{v}_{n}^{k+1}=\widehat{T}_{n}^{k}v$.

We can again apply the stochastic dominance techniques we have developed to  the convergence analysis of stochastic
process $\left\{ \|\hat{v}_{n}^{k}-v^{*}\|_{\mathcal{C}\left(\mathcal{X}\right)}\right\} _{k\geq0}$.
Similarly to that of equation \eqref{eq:defn-kappa*}, we get an upper bound
\[
\|v\|_{\mathcal{C}\left(\mathcal{X}\right)}\leq\kappa^{*}\triangleq\frac{c\, q_{\max}+\max\left\{ h\, x_{\max},\, b\, x_{\min}\right\} }{1-\alpha}
\]
for the norm of the value function of any policy for Problem (\ref{NEWSVENDOR}).
By the triangle inequality,
\[
\|\hat{v}_{n}^{k}-v^{*}\|_{\mathcal{C}\left(\mathcal{X}\right)}\leq\|\hat{v}_{n}^{k}\|_{\mathcal{C}\left(\mathcal{X}\right)}+\|v^{*}\|_{\mathcal{C}\left(\mathcal{X}\right)}\leq2\,\kappa^{*}.
\]
We can thus restrict $\|\hat{v}_{n}^{k}-v^{*}\|_{\mathcal{C}\left(\mathcal{X}\right)}$
to the state space $\left[0,\,2\,\kappa^{*}\right]$. For a fixed
granularity $\epsilon_{g}>0$, we can define $\left\{ X_{n}^{k}\right\} _{k\geq0}$
and $\left\{ Y_{n}^{k}\right\} _{k\geq0}$  as in   Section \ref{sec:randomops}. 
Our upper bound on probability follows.
\begin{proposition}
\label{prop:Hoeffding-1} For any $n\geq1$ and $\epsilon>0$
\begin{align*}
P\left\{ \|\widehat{T}_{n}v-T\, v\|\geq\epsilon\right\} \leq\, & P\left\{ \alpha\sup_{x\in\mathcal{X},\, q\in\left[0,\, q_{\max}\right]}|\mathbb{E}\left[v\left(\psi\left(x+q-D\right)\right)\right]-\frac{1}{n}\sum_{i=1}^{n}v\left(\psi\left(x+q-D_{i}\right)\right)|\geq\epsilon\right\} \\
\leq\, & 2\left\lceil \frac{9\left(L^{*}\right)^{2}q_{\max}^{2}}{\epsilon^{2}}\right\rceil e^{-2\,\left(\epsilon/3\right)^{2}n/\left(2\,\|v\|_{\mathcal{C}\left(\mathcal{X}\right)}\right)^{2}},
\end{align*}
for all $v\in\mathcal{C}\left(\mathcal{X}\right)$ with Lipschitz
constant bounded by $L^{*}$.
\end{proposition}
\proof{Proof:}
By Fact \ref{fact}, we know that

\[
\|\widehat{T}_{n}v-T\, v\|_{\mathcal{C}\left(\mathcal{X}\right)}\leq\alpha\sup_{x\in\mathcal{X},\, q\in\left[0,\, q_{\max}\right]}|\mathbb{E}\left[v\left(\psi\left(x+q-D\right)\right)\right]-\frac{1}{n}\sum_{i=1}^{n}v\left(\psi\left(x+q-D_{i}\right)\right)|.
\]
Let $\left\{ \left(x_{j},q_{j}\right)\right\} _{j=1}^{J}$ be an $\epsilon/\left(3\, L^{*}\right)-$net
for $\mathcal{X}\times\left[0,\, q_{\max}\right]$. We can choose
$J$ to be the smallest integer greater than or equal to 
\[
\frac{x_{\max}-x_{\min}}{\epsilon/\left(3\, L^{*}\right)}\times\frac{q_{\max}}{\epsilon/\left(3\, L^{*}\right)}=\frac{9\left(L^{*}\right)^{2}q_{\max}^{2}}{\epsilon^{2}}.
\]
If we have
\[
|\mathbb{E}\left[v\left(\psi\left(x_{j}+q_{j}-D\right)\right)\right]-\frac{1}{n}\sum_{i=1}^{n}v\left(\psi\left(x_{j}+q_{j}-D_{i}\right)\right)|\leq\epsilon/3
\]
for all $j=1,\ldots,J$, then
\[
|\mathbb{E}\left[v\left(\psi\left(x+q-D\right)\right)\right]-\frac{1}{n}\sum_{i=1}^{n}v\left(\psi\left(x+q-D_{i}\right)\right)|\leq\epsilon
\]
for all $\left(x,q\right)\in\mathcal{X}\times\left[0,\, q_{\max}\right]$
by Lipschitz continuity and construction of $\left\{ \left(x_{j},q_{j}\right)\right\} _{j=1}^{J}$.
Then, by Hoeffding's inequality and using union bound, we get,
\[
P\left\{ \alpha\sup_{j=1,\ldots,J}|\mathbb{E}\left[v\left(\psi\left(x_{j}+q_{j}-D\right)\right)\right]-\frac{1}{n}\sum_{i=1}^{n}v\left(\psi\left(x_{j}+q_{j}-D_{i}\right)\right)|\geq\epsilon/3\right\} \leq2\, J\, e^{-2\,\left(\epsilon/3\right)^{2}n/\left(2\,\|v\|_{\mathcal{C}\left(\mathcal{X}\right)}\right)^{2}}.
\]
\hfill \Halmos \endproof
As before, we use the preceding complexity estimate to determine $p_{n}$
for the family $\left\{ Y_{n}^{k}\right\} _{k\geq0}$. The remainder
of our stochastic dominance argument is exactly the same.

\section{Numerical Experiments}\label{sec:experimental}

\begin{figure}
\FIGURE
{\includegraphics[width=\textwidth]{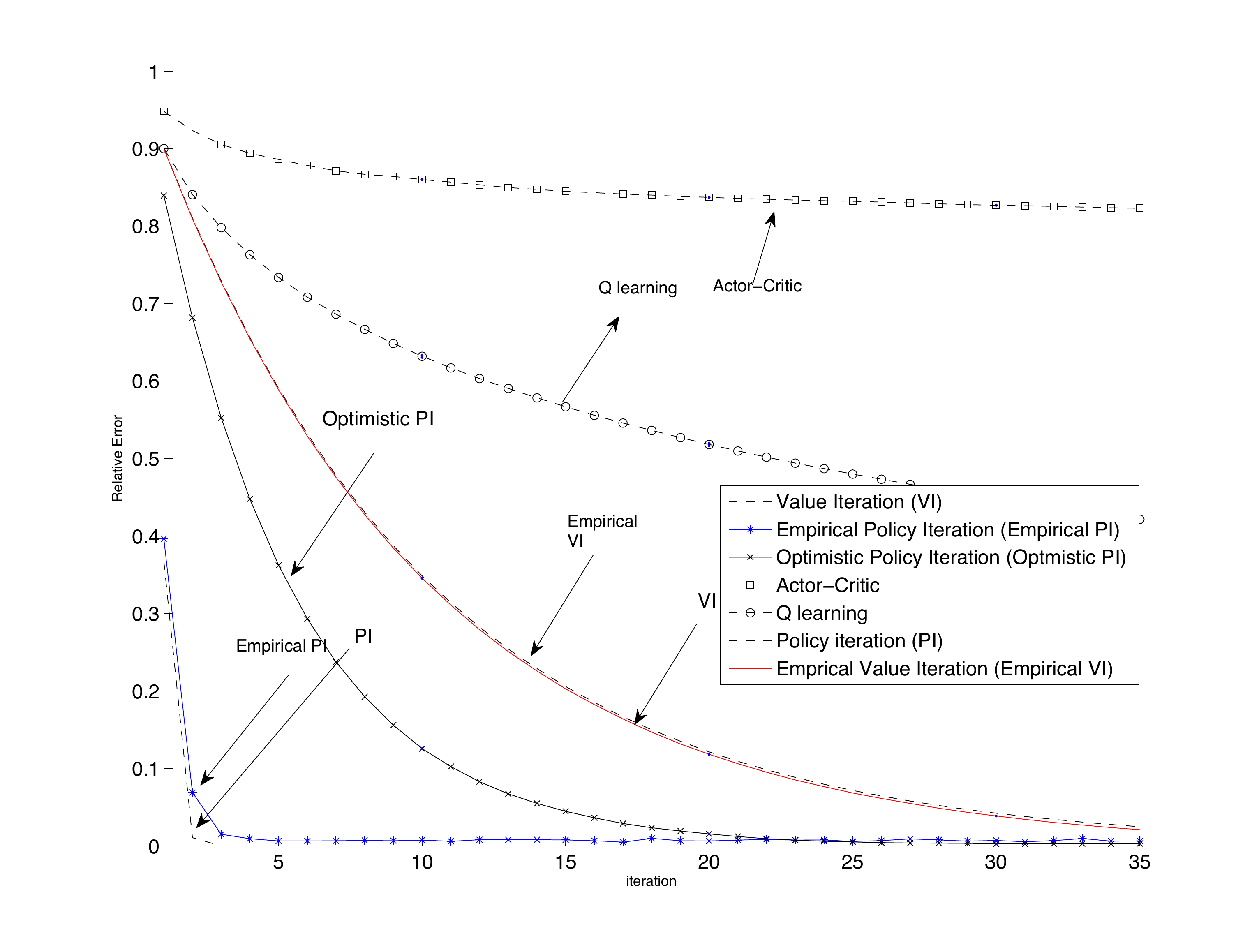}}
{Numerical performance comparison of empirical value and policy iteration with actor-critic, Q-learning and optimistic policy iteration. Number of samples taken: $n=10, q=10$ (Refer EVI and EPI algorithms in Section 3 for the definition of $n$ and $q$). \label{fig:EDP-comparison}}
{}
\end{figure}

We now provide a numerical comparison of EDP methods with other methods for approximate dynamic programming via simulation.
Figure \ref{fig:EDP-comparison} shows relative error ($||\hat{v}_n^k-v^*||$) of the Actor-Critic algorithm, Q-Learning algorithm, Optimistic Policy Iteration (OPI), Empirical Value Iteration (EVI) and Empirical Policy Iteration (EPI). It also shows relative error for exact Value iteration (VI). The problem considered was a generic 10 state and 5 action space MDP with infinite-horizon discounted cost. From the figure, we see that EVI and EPI significantly outperform Actor-Critic algorithm (which converges very slowly) and Q-Learning. Optimistic policy iteration performs better than EVI since policy iteration-based algorithms are known to converge faster, but EPI outperforms OPI as well. 

The experiments were peformed on a generic laptop with Intel Core i7 processor, 4GM RAM, on a 64-bit Windows 7 operating system via Matlab R2009b environment. 

These preliminary numerical results seem to suggest that EDP methods outperform other ADP methods numerically, and hold good promise. More definitive conclusions about their numerical performance requires further work. We would also like to point that EDP methods would very easily be parallelizable, and hence, they could potentially be useful for a wider variety of problem settings.

\section{Conclusions}\label{sec:conclusions}

In this paper, we have introduced a new class of algorithms for approximate dynamic programming. The idea is actually not novel, and quite simple and natural: just replace the expectation in the Bellman operator with an empirical estimate (or a sample average approximation, as it is often called.) The difficulty, however, is that it makes the Bellman operator a random operator. This makes its convergence analysis very challenging since (infinite horizon) dynamic programming theory is based on looking at the fixed points of the Bellman operator. However, the extant notions of `probabilistic' fixed points for random operators are not relevant since they are akin to classical fixed points of deterministic monotone operators when $\omega$ is fixed. We introduce two notions of probabilistic fixed points - strong and weak. Furthermore, we show that these asymptotically coincide with the classical fixed point of the related deterministic operator, This is reassuring as it suggests that approximations to our probabilistic fixed points (obtained by finitely many iterations of the empirical Bellman operator)  are going to be approximations to the classical fixed point of the Bellman operator as well.

In developing this theory, we also developed a mathematical framework based on stochastic dominance for convergence analysis of random operators. While our immediate goal was analysis of iteration of the empirical Bellman operator in empirical dynamic programming, the framework is likely of broader use, possibly after further development.

We have then shown that many variations and extensions of classical dynamic programming, work for empirical dynamic programming as well. In particular, empirical dynamic programming can be done asynchronously just as classical DP can be. Moreover, a zero-sum stochastic game can be solved by a minimax empirical dynamic program. We also apply the EDP method to the dynamic newsvendor problem which has continuous state and action spaces, which demonstrates the potential of EDP to solve problems over more general state and action spaces. 

We have done some preliminary experimental performance analysis of EVI and EPI, and compared it to similar methods. Our numerical simulations suggest that EDP algorithms converge faster than stochastic approximation-based actor-critic, Q-learning and optimistic policy iteration algorithms. However, these results are only suggestive, we do not claim definitive performance improvement in practice over other algorithms. This requires an extensive and careful numerical investigation of all such algorithms. 

We do note that EDP methods, unlike stochastic approximation methods, do not require any recurrence property to hold. In that sense, they are more universal. On the other hand, EDP algorithms would inherit some of the `curse of dimensionality' problems associated with exact dynamic programming. Overcoming that challenge requires additional ideas, and is potentially a fruitful direction for future research. Some other directions of research are extending the EDP algorithms to the infinite-horizon average cost case, and to the partially-observed case. We will take up these issues in the future.

\paragraph*{Acknowledgements}
The authors would like to thank Ugur Akyol (USC) for running the numerical experiments for this research. The authors would also like to thank Suvrajeet Sen (USC) and Vivek Borkar (IIT Bombay) for initial feedback on this work.

\newpage
\bibliographystyle{ormsv080}
\bibliography{References1}

\newpage
\begin{APPENDIX}{}

\section{Proofs of Various Lemmas, Propositions and Theorems} 
\subsection{Proof of Fact \ref{fact}}
\label{pf:fact}
\proof{Proof:}
To verify part (i), note

\begin{align*}
\inf_{x\in X}f_{1}\left(x\right)=\, & \inf_{x\in X}\left\{ f_{1}\left(x\right)+f_{2}\left(x\right)-f_{2}\left(x\right)\right\} \\
\leq\, & \inf_{x\in X}\left\{ f_{2}\left(x\right)+|f_{1}\left(x\right)-f_{2}\left(x\right)|\right\} \\
\leq\, & \inf_{x\in X}\left\{ f_{2}\left(x\right)+\sup_{y\in Y}|f_{1}\left(y\right)-f_{2}\left(y\right)|\right\} \\
\leq\, & \inf_{x\in X}f_{2}\left(x\right)+\sup_{y\in Y}|f_{1}\left(y\right)-f_{2}\left(y\right)|,
\end{align*}
giving
\[
\inf_{x\in X}f_{1}\left(x\right)-\inf_{x\in X}f_{2}\left(x\right)\leq\sup_{x\in X}|f_{1}\left(x\right)-f_{2}\left(x\right)|.
\]
By the same reasoning,

\[
\inf_{x\in X}f_{2}\left(x\right)-\inf_{x\in X}f_{1}\left(x\right)\leq\sup_{x\in X}|f_{1}\left(x\right)-f_{2}\left(x\right)|,
\]
and the preceding two inequalities yield the desired result. Part
(ii) follows similarly.
\hfill \Halmos \endproof

\subsection{Proof of Theorem \ref{thm:dominance}}
\label{pf:thm:dominance}

We first prove the following lemmas. 

%

\begin{lemma}
\label{lem:ST_monotonicity} $\left[Y_{n}^{k+1}\,\vert\, Y_{n}^{k}=\theta\right]$
is stochastically increasing in $\theta$ for all $k\geq0$, i.e.\\
$\left[Y_{n}^{k+1}\,\vert\, Y_{n}^{k}=\theta\right]\leq_{st}\left[Y_{n}^{k+1}\,\vert\, Y_{n}^{k}=\theta'\right]$
for all $\theta\leq\theta'$.\end{lemma}
\proof{Proof:}
We see that
\[
\mbox{Pr}\left\{ Y_{n}^{k+1}\geq\eta\,\vert\, Y_{n}^{k}=\theta\right\} 
\]
is increasing in $\theta$ by construction of $\left\{ Y_{n}^{k}\right\} _{k\geq0}$.
If $\theta>\eta$, then $\mbox{Pr}\left\{ Y_{n}^{k+1}\geq\eta\,\vert\, Y_{n}^{k}=\theta\right\} =1$
since $Y_{n}^{k+1}\geq\theta-1$ almost surely; if $\theta\leq\eta$,
then $\mbox{Pr}\left\{ Y_{n}^{k+1}\geq\eta\,\vert\, Y_{n}^{k}=\theta\right\} =1-p_{n}$
since the only way $Y_{n}^{k+1}$ will remain larger than $\eta$
is if $Y_{n}^{k+1}=N^{*}$.
\hfill \Halmos \endproof

\begin{lemma}
\label{lem:ST_XtYt} 
$\left[X_{n}^{k+1}\,\vert\, X_{n}^{k}=\theta,\mathcal{F}^{k}\right]\leq_{st}\left[Y_{n}^{k+1} \vert\, Y_{n}^{k}=\theta\right]$
for all $\theta$ and all $\mathcal{F}^{k}$ for all $k\geq0$.\end{lemma}
\proof{Proof:}
Follows from construction of $\left\{ Y_{n}^{k}\right\} _{k\geq0}$.
For any history $\mathcal{F}^{k}$,
\begin{align*}
\mathcal{P}\left\{ X_{n}^{k+1}\geq\theta-1\,\vert\, X_{n}^{k}=\theta,\,\mathcal{F}^{k}\right\} \leq \mathcal{Q}\left\{ Y_{n}^{k+1}\geq\theta-1\,\vert\, Y_{n}^{k}=\theta\right\} =1.
\end{align*}
Now, 
\begin{align*}
\mathcal{P}\left\{ X_{n}^{k+1}=N^{*}\,\vert\, X_{n}^{k}=\theta,\,\mathcal{F}^{k}\right\}  &\leq \mathcal{P}\left\{ X_{n}^{k+1} > \theta - 1 \,\vert\, X_{n}^{k}=\theta,\,\mathcal{F}^{k}\right\} \\
&= 1 - P(X_{n}^{k+1}  \leq \theta -1 | X^{k}_{n} = \theta) \\
&\leq 1 - p_{n}
\end{align*}
the last inequality  follows because $p_{n}$ is the worst case probability for a one-step improvement in the Markov chain $\{X^{k}_{n}\}_{k \geq 0}$ 
\hfill \Halmos \endproof

\noindent {\bf Proof of Theorem \ref{thm:dominance}}\\
\proof{Proof:}
Trivially, $X_{n}^{0}\leq_{st}Y_{n}^{0}$ since $X_{n}^{0}=Y_{n}^{0}$.
Next, we see that $X_{n}^{1}\leq_{st}Y_{n}^{1}$ by previous lemma. 
We prove the general case by induction. Suppose $X_{n}^{k}\leq_{st}Y_{n}^{k}$
for $k\geq1$, and for this proof define the random variable
\[
\mathfrak{Y}\left(\theta\right)=\begin{cases}
\max\left\{ \theta-1,\,\eta^{*}\right\} , & \mbox{w.p. }p_{n},\\
N^{*}, & \mbox{w.p. }1-p_{n},
\end{cases}
\]
as a function of $\theta$. We see that $Y_{n}^{k+1}$ has the same
distribution as
\[
\left[\mathfrak{Y}\left(\Theta\right) | \Theta=Y_{n}^{k}\right]
\]
by definition. Since $\mathfrak{Y}\left(\theta\right)$ are stochastically
increasing, we see that
\[
\left[\mathfrak{Y}\left(\Theta\right) | \Theta=Y_{n}^{k}\right]\geq_{st}\,\left[\mathfrak{Y}\left(\Theta\right) | \Theta=X_{n}^{k}\right]
\]
by \cite[Theorem 1.A.6]{shaked2007stochastic} and our induction hypothesis.
Now,
\[
\left[\mathfrak{Y}\left(\Theta\right) | \Theta=X_{n}^{k}\right]\geq_{st}\left[X_{n}^{k+1} | X_{n}^{k},\,\mathcal{F}^{k}\right]
\]
by \cite[Theorem 1.A.3(d)]{shaked2007stochastic} for all histories $\mathcal{F}^{k}$.
It follows that $Y_{n}^{k+1}\geq_{st}X_{n}^{k+1}$ by transitivity.
\hfill \Halmos \endproof

\subsection{Proof of Lemma \ref{lem:steadystate}}
\label{pf:lem:steadystate}

\proof{Proof:}
The stationary probabilities $\left\{ \mu_{n}\left(i\right)\right\} _{i=\eta^{*}}^{N^{*}}$
satisfy the equations:
\begin{align*}
\mu_{n}\left(\eta^{*}\right)=\, & p_{n}\mu_{n}\left(\eta^{*}\right)+p_{n}\mu_{n}\left(\eta^{*}+1\right),\\
\mu_{n}\left(i\right)=\, & p_{n}\mu_{n}\left(i+1\right), & \forall i=\eta^{*}+1,\ldots,N^{*}-1,\\
\mu_{n}\left(N^{*}\right)=\, & \left(1-p_{n}\right)\sum_{i=\eta^{*}}^{N^{*}}\mu_{n}\left(i\right),\\
\sum_{i=\eta^{*}}^{N^{*}}\mu_{n}\left(i\right)=\, & 1.
\end{align*}
We then see that

\[
\mu_{n}\left(i\right)=p_{n}^{\left(N^{*}-i\right)}\mu_{n}\left(N^{*}\right),\,\forall i=\eta^{*}+1,\ldots,N^{*}-1,
\]
and
\[
\mu_{n}\left(\eta^{*}\right)=\frac{p_{n}}{1-p_{n}}\mu_{n}\left(\eta^{*}+1\right)=\frac{p_{n}^{N^{*}-\eta^{*}}}{1-p_{n}}\mu\left(N^{*}\right).
\]
We can solve for $\mu_{n}\left(N^{*}\right)$ using $\sum_{i=\eta^{*}}^{N^{*}}\mu_{n}\left(i\right)=1$,
\begin{align*}
1=\, & \sum_{i=\eta^{*}}^{N^{*}}\mu_{n}\left(i\right)\\
=\, & \frac{p_{n}^{N^{*}-\eta^{*}}}{1-p_{n}}\mu_{n}\left(N^{*}\right)+\sum_{i=\eta^{*}+1}^{N^{*}}p_{n}^{N^{*}-i}\mu_{n}\left(N^{*}\right)\\
=\, & \left[\frac{p_{n}^{N^{*}-\eta^{*}}}{1-p_{n}}+\sum_{i=\eta^{*}+1}^{N^{*}}p_{n}^{N^{*}-i}\right]\mu_{n}\left(N^{*}\right)\\
=\, & \left[\frac{p_{n}^{N^{*}-\eta^{*}}}{1-p_{n}}+\frac{p_{n}-p_{n}^{N^{*}-\eta^{*}}}{1-p_{n}}\right]\mu_{n}\left(N^{*}\right),\\
=\, & \frac{p_{n}}{1-p_{n}}\mu_{n}\left(N^{*}\right),
\end{align*}
based on the fact that
\[
\sum_{i=\eta^{*}+1}^{N^{*}}p_{n}^{\left(N^{*}-i\right)}=\sum_{i=0}^{N^{*}-\eta^{*}-1}p_{n}^{i}=\frac{1-p_{n}^{\left(N^{*}-\eta^{*}\right)}}{1-p_{n}}.
\]
We conclude
\[
\mu\left(N^{*}\right)=\frac{1-p_{n}}{p_{n}},
\]
and thus
\[
\mu_{n}\left(i\right)=\left(1-p_{n}\right)p_{n}^{\left(N^{*}-i-1\right)},\,\forall i=\eta^{*}+1,\ldots,N^{*}-1,
\]
and
\[
\mu_{n}\left(\eta^{*}\right)=p_{n}^{N^{*}-\eta^{*}-1}. 
 \]
 \hfill \hfill \Halmos 
\endproof

\subsection{Proof of Proposition \ref{lem:strongpfp}}
\label{pf:lem:strongpfp}
\proof{Proof:}
(i) First observe that 
\[
\lim_{n\rightarrow\infty}\widehat{T}_{n}\left(\omega\right)v^{*}=T\, v^{*},
\]
by Assumption \ref{ass:fixedpoint}. It follows that $\widehat{T}_{n}\left(\omega\right)v^{*}$
converges to $v^{*}=T\, v^{*}$ as $n\rightarrow\infty$, $P-$almost
surely. Almost sure convergence implies convergence in probability.\\
(ii) Let $\hat{v}$ be a strong probabilistic fixed point. Then,
\[P(\|T \hat{v} - \hat{v}\| \geq \epsilon)  \leq P(\|T \hat{v} - \widehat{T}_{n} \hat{v}\| \geq \epsilon/2) + P(\|\widehat{T}_{n} \hat{v} - \hat{v}\| \geq \epsilon/2 ) \]
First term on the RHS can be made arbitrarily small by Assumption \ref{ass:fixedpoint}. Second term on RHS can  also be made arbitrarily small   by the definition of strong probabilistic fixed point. So, for sufficiently large $n$, we get $P(\|T \hat{v} - \hat{v}\| \geq \epsilon) < 1$. Since the event in the LHS is deterministic, we get $\|T \hat{v} - \hat{v}\|=0$. Hence, $\hat{v} =  v^{*}$. 

\hfill \Halmos

\endproof

\subsection{Proof of Proposition \ref{prop:wkpfp1}}
\label{pf:prop:wkpfp1}
\proof{Proof:}
(i) This statement is proved in Theorem \ref{thm:asymptotic}.

(ii) Fix the initial seed $v\in\mathbb{R}^{|\mathbb{S}|}$. For a
contradiction, suppose $\hat{v}$ is not a fixed point of $T$ so
that $\|v^{*}-\hat{v}\|=\epsilon'>0$ (we use here the fact that $v^{*}$
is unique). Now
\[
\|\hat{v}-v^{*}\|=\epsilon'\leq\|\widehat{T}_{n}^{k}v-\hat{v}\|+\|\widehat{T}_{n}^{k}v-v^{*}\|
\]
for any $n$ and $k$ by the triangle inequality. For clarity, this
inequality holds in the almost sure sense:
\[
\mathcal{P}\left(\epsilon'\leq\|\widehat{T}_{n}^{k}v-\hat{v}\|+\|\widehat{T}_{n}^{k}v-v^{*}\|\right)=1
\]
for all $n$ and $k$.

We already know that
\[
\lim_{n\rightarrow\infty}\limsup_{k\rightarrow\infty}\mathcal{P}\left(\|\widehat{T}_{n}^{k}v-v^{*}\|>\epsilon'/3\right)=0
\]
by Theorem \ref{thm:asymptotic} and
\[
\lim_{n\rightarrow\infty}\limsup_{k\rightarrow\infty}\mathcal{P}\left(\|\widehat{T}_{n}^{k}v-\hat{v}\|>\epsilon'/3\right)=0
\]
by assumption. Now
\[
\mathcal{P}\left(\max\left\{ \|\widehat{T}_{n}^{k}v-\hat{v}\|,\,\|\widehat{T}_{n}^{k}v-v^{*}\|\right\} >\epsilon'/3\right)\leq\mathcal{P}\left(\|\widehat{T}_{n}^{k}v-v^{*}\|>\epsilon'/3\right)+\mathcal{P}\left(\|\widehat{T}_{n}^{k}v-\hat{v}\|>\epsilon'/3\right),
\]
so
\[
\lim_{n\rightarrow\infty}\limsup_{k\rightarrow\infty}\mathcal{P}\left(\max\left\{ \|\widehat{T}_{n}^{k}v-\hat{v}\|,\,\|\widehat{T}_{n}^{k}v-v^{*}\|\right\} >\epsilon'/3\right)=0.
\]
However, $\epsilon'\leq\|\widehat{T}_{n}^{k}v-\hat{v}\|+\|\widehat{T}_{n}^{k}v-v^{*}\|$
almost surely so at least one of $\|\widehat{T}_{n}^{k}v-\hat{v}\|$
or $\|\widehat{T}_{n}^{k}v-v^{*}\|$ must be greater than $\epsilon'/3$
for all large $k$.
\hfill \Halmos \endproof

\subsection{Proof of Proposition \ref{prop:empB-ass}}
\label{pf:prop:empB-ass}
\proof{Proof:}
(i) Assumption \ref{ass:fixedpoint} : \\
Certainly, 
\[
\|T_{n}\left(\omega\right)\left(v\right)-T\, v\|\leq\alpha\max_{\left(s,a\right)\in\mathbb{K}}|\frac{1}{n}\sum_{i=1}^{n}v\left(\psi\left(s,a,\xi_{i}\right)\right)-\mathbb{E}\left[v\left(\psi\left(s,a,\xi\right)\right)\right]|
\]
using Fact \eqref{fact}. We know that for any fixed $\left(s,a\right)\in\mathbb{K}$,

\[
|\frac{1}{n}\sum_{i=1}^{n}v\left(\psi\left(s,a,\xi_{i}\right)\right)-\mathbb{E}\left[v\left(\psi\left(s,a,\xi\right)\right)\right]|\rightarrow0,
\]
as $n\rightarrow\infty$ by the Strong Law of Large Numbers (the random
variable $v\left(\psi\left(s,a,\xi\right)\right)$ has finite expectation
because it is essentially bounded). Recall that $\mathbb{K}$ is finite
to see that the right hand side of the above inequality converges
to zero as $n\rightarrow\infty$.\\
(ii) Assumption \ref{ass:bddness} :\\
We define the constant
\[
\kappa^{*}\triangleq\frac{\max_{\left(s,a\right)\in\mathbb{K}}|c\left(s,a\right)|}{1-\alpha}.
\]
Then it can be easily verified that the value of  any policy $\pi$ $v^{\pi} \leq \kappa^{*}$. Then $v^{*} \leq \kappa^{*}$ and without loss of generality we can restrict $\hat{v}^{k}_{n}$ to the set $\overline{B_{2\kappa^{*}}\left(0\right)}$.\\
(iii) Assumption \ref{ass:rate}:\\
This is the well known contraction property of the Bellman operator.\\
(iv) Assumption \ref{ass:pn}:\\
Using Fact \ref{fact}, for any fixed $s \in \mathbb{S}$, 
\begin{align*}
|\widehat{T}_{n}v(s) &- T v(s)| \leq \max_{a \in A(s)} \left|\frac{\alpha}{n} \sum^{n}_{i=1}v(\psi(s, a, \xi_{i}) - \alpha \mathbb{E}[v(\psi(s, a, \xi)] \right|
\end{align*}
and hence,
\begin{align*}
&P\left\{ \|\widehat{T}_{n}v -T\, v\|\geq\epsilon\right\}  \leq P\left\{\max_{(s,a) \in \mathbb{K}} \left|\frac{\alpha}{n} \sum^{n}_{i=1}v(\psi(s, a, \xi_{i}) - \alpha \mathbb{E}[v(\psi(s, a, \xi)] \right| \geq \epsilon  \right\}.
\end{align*}

For any fixed $\left(s,a\right)\in\mathbb{K}$,
\begin{align*}
P\left\{ \left |\frac{\alpha}{n}\sum_{i=1}^{n}v\left(\psi\left(s,a,\xi_{i}\right)\right)-\alpha \mathbb{E}\left[v\left(\psi\left(s,a,\xi\right)\right)\right]\right|\geq\epsilon\right\} &\leq 2\, e^{-2\,(\epsilon/\alpha)^{2}n/\left(v_{\max}-v_{\min}\right)^{2}}\leq2\, e^{-2\,(\epsilon/\alpha)^{2}n/\left(2\,\|v\|\right)^{2}}  \\
& \leq 2\, e^{-2\,(\epsilon/\alpha)^{2}n/\left(2\,\kappa^{*}\right)^{2}}  
\end{align*}
by  Hoeffding's inequality. Then, using the union bound, we have
\[P\left\{ \|\widehat{T}_{n}v -T\, v\|\geq\epsilon\right\}  \leq  2\,|\mathbb{K}|\,e^{-2\,(\epsilon/\alpha)^{2}n/\left(2\,\kappa^{*}\right)^{2}}.  \]
By taking complements of the above event we get the desired result. 
\hfill \Halmos 
\endproof

%

\subsection{Lemma \ref{lem:eigen} }
\label{pf:lem:eigen}

%
%
%

\begin{lemma}
\label{lem:eigen} 
For any fixed $n \geq 1$, the eigenvalues of the transition matrix $\mathfrak{Q}$ of the Markov chain $Y^{k}_{n}$ are $0$ (with algebraic multiplicity $N^{*}-\eta^{*}-1$) and 1.
\end{lemma}

\proof{Proof:}
In general,
the transition matrix $\mathfrak{Q}_{n}\in\mathbb{R}^{\left(N^{*}-\eta^{*}+1\right)\times\left(N^{*}-\eta^{*}+1\right)}$
of $\left\{ Y_{n}^{k}\right\} _{k\geq0}$ looks like 
\[
\mathfrak{Q}_{n}=\left[\begin{array}{cccccc}
p_{n} & 0 & \cdots & \cdots & 0 & \left(1-p_{n}\right)\\
p_{n} & 0 & \cdots & \cdots & 0 & \left(1-p_{n}\right)\\
0 & p_{n} & 0 & \cdots & 0 & \left(1-p_{n}\right)\\
\vdots & \vdots & \vdots & \ddots & \vdots & \vdots\\
0 & 0 & \cdots & \cdots & 0 & \left(1-p_{n}\right)\\
0 & 0 & \cdots & \cdots & p_{n} & \left(1-p_{n}\right)
\end{array}\right].
\]

To compute the eigenvalues of $\mathfrak{Q}_{n}$, we want to solve
$\mathfrak{Q}_{n}x=\lambda\, x$ for some $x\ne0$ and $\lambda\in\mathbb{R}$.
For $x=\left(x_{1},\, x_{2},\ldots,\, x_{N^{*}-\eta^{*}+1}\right)\in\mathbb{R}^{N^{*}-\eta^{*}+1}$,
\[
\mathfrak{Q}_{n}x=\left(\begin{array}{c}
p_{n}x_{1}+\left(1-p_{n}\right)x_{N^{*}-\eta^{*}+1}\\
p_{n}x_{1}+\left(1-p_{n}\right)x_{N^{*}-\eta^{*}+1}\\
p_{n}x_{2}+\left(1-p_{n}\right)x_{N^{*}-\eta^{*}+1}\\
\vdots\\
p_{n}x_{N^{*}-\eta^{*}-1}+\left(1-p_{n}\right)x_{N^{*}-\eta^{*}+1}\\
p_{n}x_{N^{*}-\eta^{*}}+\left(1-p_{n}\right)x_{N^{*}-\eta^{*}+1}
\end{array}\right).
\]
Now, suppose $\lambda\ne0$ and $\mathfrak{Q}\, x=\lambda\, x$ for
some $x\ne0$. By the explicit computation of $\mathfrak{Q}\, x$
above,
\[
\left[\mathfrak{Q}\, x\right]_{1}=p\, x_{1}+\left(1-p\right)x_{N^{*}-\eta^{*}+1}=\lambda\, x_{1}
\]
and
\[
\left[\mathfrak{Q}\, x\right]_{2}=p\, x_{1}+\left(1-p\right)x_{N^{*}-\eta^{*}+1}=\lambda\, x_{2},
\]
so it must be that $x_{1}=x_{2}$. However, then
\[
\left[\mathfrak{Q}\, x\right]_{2}=p\, x_{1}+\left(1-p\right)x_{N^{*}-\eta^{*}+1}=p\, x_{2}+\left(1-p\right)x_{N^{*}-\eta^{*}+1}=\left[\mathfrak{Q}\, x\right]_{3},
\]
and thus $x_{2}=x_{3}$. Continuing this reasoning inductively shows
that $x_{1}=x_{2}=\cdots=x_{N^{*}-\eta^{*}+1}$ for any eigenvector
$x$ of $\mathfrak{Q}$. Thus, it must be true that $\lambda=1$.
\hfill \Halmos \endproof

%
%
%

\subsection{Proof of Lemma \ref{lem:asyn-Tx}}
\label{pf:lem:asyn-Tx}

\proof{Proof:}
(i) Suppose $v\leq v'$. It is automatic that $\left[T_{x}v\right]\left(s\right)=v\left(s\right)\leq v'\left(s\right)=\left[T_{x}v'\right]\left(s\right)$
for all $s\ne x$. For state $s=x$,
\[
c\left(s,a\right)+\alpha\,\mathbb{E}\left[v\left(s'\right)|s,a\right]\leq c\left(s,a\right)+\alpha\,\mathbb{E}\left[v'\left(s'\right)|s,a\right]
\]
for all $\left(s,a\right)\in\mathbb{K}$, so
\[
\min_{a\in A\left(s\right)}\left\{ c\left(s,a\right)+\alpha\,\mathbb{E}\left[v\left(s'\right)|s,a\right]\right\} \leq\min_{a\in A\left(s\right)}\left\{ c\left(s,a\right)+\alpha\,\mathbb{E}\left[v'\left(s'\right)|s,a\right]\right\} ,
\]
and thus $\left[T_{x}v\right]\left(s\right)\leq\left[T_{x}v'\right]\left(s\right)$.

(ii) We see that
\[
\min_{a\in A\left(s\right)}\left\{ c\left(s,a\right)+\alpha\,\mathbb{E}\left[v'\left(s'\right)+\eta|s,a\right]\right\} =\min_{a\in A\left(s\right)}\left\{ c\left(s,a\right)+\alpha\,\mathbb{E}\left[v'\left(s'\right)|s,a\right]\right\} +\alpha\,\eta
\]
for state $x$ and all other states are not changed.
\hfill \Halmos \endproof

\subsection{Proof of Proposition \ref{prop:Asynchronous_ErrorBound}}
\label{pf:prop:Asynchronous_ErrorBound}
\proof{Proof:}
Starting with $v^{0}$,
\[
T_{x_{0}}v^{0}-\eta\,1\leq T_{x_{0}}v^{0}+\varepsilon_{0}\leq T_{x_{0}}v^{0}+\eta\,1,
\]
which gives

\[
T_{x_{0}}v^{0}-\eta\,1\leq\tilde{v}^{0}\leq T_{x_{0}}v^{0}+\eta\,1,
\]
and 
\[
v^{1}-\eta\,1\leq\tilde{v}^{1}\leq v^{1}+\eta\,1.
\]
By monotonicity of $T_{x_{1}}$,
\[
T_{x_{1}}\left[v^{1}-\eta\,1\right]\leq T_{x_{1}}\tilde{v}^{1}\leq T_{x_{1}}\left[v^{1}+\eta\,1\right],
\]
and by our assumptions on the noise,
\[
T_{x_{1}}\left[v^{1}-\eta\,1\right]-\eta\,1\leq\tilde{v}^{2}\leq T_{x_{1}}\left[v^{1}+\eta\,1\right]+\eta\,1.
\]
Now
\[
T_{x_{1}}\left[v-\eta\,1\right]=T_{x_{1}}v-\alpha\,\eta\, e_{x},
\]
thus
\[
v^{2}-\eta\,1-\alpha\,\eta\,1\leq\tilde{v}^{2}\leq v^{2}+\eta\,1+\alpha\,\eta\,1.
\]
Similarly,
\[
v^{3}-\alpha\left(\eta\,1-\alpha\,\eta\,1\right)-\eta\,1\leq\tilde{v}^{3}\leq v^{3}+\alpha\left(\eta\,1+\alpha\,\eta\,1\right)+\eta\,1,
\]
and the general case follows.
\hfill \Halmos \endproof

\subsection{Proof of Lemma \ref{lem:ShapleyT}}
\label{pf:lem:ShapleyT}
\proof{Proof:}
Using Fact \ref{fact} twice, compute
\begin{align*}
 |\left[T\, v\right](s)-\left[T\, v'\right](s)| =\, & |\min_{a\in A\left(s\right)}\max_{b\in B\left(s\right)}\left\{ r\left(s,a,b\right)+\alpha\,\mathbb{E}\left[v\left(\tilde{s}\right)|s,a,b\right]\right\} \\
&\hspace{2cm}-\min_{a\in A\left(s\right)}\max_{b\in B\left(s\right)}\left\{ r\left(s,a,b\right)+\alpha\,\mathbb{E}\left[v'\left(\tilde{s}\right)|s,a,b\right]\right\} |\\
\leq\, & \max_{a\in A}|\max_{b\in B\left(s\right)}\left\{ r\left(s,a,b\right)+\alpha\,\mathbb{E}\left[v\left(\tilde{s}\right)|s,a,b\right]\right\} -\max_{b\in B\left(s\right)}\left\{ r\left(s,a,b\right)+\alpha\,\mathbb{E}\left[v'\left(\tilde{s}\right)|s,a,b\right]\right\} |\\
\leq\, & \alpha\,\max_{a\in A}\max_{b\in B\left(s\right)}|\mathbb{E}\left[v\left(\tilde{s}\right)-v'\left(\tilde{s}\right)|s,a,b\right]|\\
\leq\, & \alpha\,\max_{a\in A}\max_{b\in B}\mathbb{E}\left[|v\left(\tilde{s}\right)-v'\left(\tilde{s}\right)||s,a,b\right]\\
\leq\, & \alpha\,\|v-v'\|.
\end{align*}

\hfill \Halmos \endproof

\end{APPENDIX}

\end{document}